\documentclass[a4paper,11pt]{amsart}
\usepackage{url}
\usepackage[colorlinks=false]{hyperref}
\usepackage{amssymb,amscd,amsmath,amsthm,amsfonts,graphicx,color}
\usepackage[all]{xy}
\usepackage{rotating}

\newcommand\boxb[1]{\square_b}

\numberwithin{equation}{section}

\newcommand\paperbody%
        {}

\newtheorem{lemma}{Lemma}

\newtheorem{theorem}{Theorem}
\newtheorem{non-theorem}{Non-Theorem}

\theoremstyle{remark}
\newtheorem{definition}{Definition}

\newcommand{\ba}{\begin{array}}
\newcommand{\ea}{\end{array}}

\begin{document}

\title[On the space of evolutionary operators and folding dynamics]
{A note on the space of evolutionary operators in population genetics and folding dynamics}

\author{Sergio Luki\'c}
\address{
\newline Corresponding address: Simons Center for Systems Biology, School of Natural Sciences, Institute for Advanced Study, Princeton NJ 08540, USA.
}
\email{lukic@ias.edu}

\begin{abstract} 
Discrete dynamical systems defined by the iteration of a polynomial map of the unit simplex to itself appear in the context of
population genetic systems evolving under mutation, recombination and weak selection. 
Although exceptional progress has been made in finding particular solutions to these systems, our knowledge 
of the general properties of the space of all possible dynamical systems of this kind is
still limited. We prove that the space of bounded-degree polynomial
maps of the unit simplex to itself is a compact and convex subset of a Euclidean space.
We provide an explicit characterization of such a space and of its boundary. 
A special class of maps in the boundary, the folding maps, which generalize the logistic map for any dimension and degree are 
defined and constructed.
Finally, we use numerical methods to study the ergodic and mixing properties of maps in the neighborhood of several of these
folding maps.
\end{abstract}


\maketitle

\setcounter{tocdepth}{3}

\paperbody

\section{Introduction}

Growing access to genome data generated by next-generation sequencing technologies is 
revitalizing the field of mathematical population genetics. In particular, the problem
of determining the simplest historical scenarios that have likely shaped patterns of genetic variation 
is motivating the development of new mathematical and computational tools that can get the most
out of these large datasets.
For instance, evolutionary biologists are often interested in ascertaining the genealogical 
relationships that underlie a sample of homologous DNA sequences.
In this case, one can model the background demographic history of 
the populations from which the samples were taken (e.g. divergence between populations, migratory flows between populations, 
changes in the population sizes, etc.), together with the effects of mutation, recombination and natural selection.
Using these kinds of models, investigators can infer the simplest scenarios that best fit the observed patterns of genetic variation and the most
probable genealogies that underlie a given set of DNA sequences \cite{SH13}.

Some of the mathematical methods that are used in these applications consist of:
\begin{itemize} 
\item[i)] Computations of the likelihood function by means of integrals defined over the space of all possible genealogies 
associated with the sample \cite{Fel88}.
\item[ii)] Sequentially Markovian approximations. These methods capture the complexities of the genealogical relationships
by means of Hidden Markov Models defined as stochastic processes along the DNA sequences \cite{MC05}.
\item[iii)] Computations of Allele Frequency Spectra (AFS) by means of random and/or deterministic dynamical systems defined 
on the space of population frequencies \cite{LH12}. In this case, rather than working directly with DNA sequence data, one uses
summary statistics of the data defined by the frequencies of the alternative forms of the genes in the sample.
\end{itemize}
Recent progress in the development of AFS-based techniques in addition to their computational efficiency and mathematical elegance
make this approach particularly appealing.
 
Allele frequency spectra summarize the distribution of counts of alternative forms of genes (alleles) 
that are observed in finite samples of DNA sequences.
This distribution is a function of the population frequencies of
the corresponding gene variants, or alleles. 
It is on this space of population allele frequencies in which the dynamical systems that model the evolutionary 
process is defined.
The space can be represented as the direct product of unit simplices. The dimensionality of the simplices
depends on the number of alleles associated with the genes, while the number of simplices 
in the direct product corresponds to the number of populations. For simplicity, here we consider one single population. 
If the number of variants of a gene in the population is finite and equal to $n+1$,
we can enumerate each allele by an integer number $i$ and denote by $x_i$
the proportion of individuals that share this particular form of the gene. The genetic composition of the population $\{ x_i\}_{i=1}^{n+1}$ 
is then described by a point in the $n$-dimensional unit simplex.
In the absence of genetic drift\footnote{Genetic drift denotes the sampling noise that 
arises in populations with finite number of individuals after sampling the DNA sequences that constitute each new generation of individuals.}, 
the dynamics of the allele frequencies is deterministic and very often 
can be described by means of a polynomial map from this simplex into itself. 
For instance, the dynamics of a population evolving under recurrent mutations 
are described by means of finite Markov chains.
Assuming that the population mates at random, the corresponding time evolution is simply given by a linear map
\begin{equation}
\label{blaone}
x_i(t+1)= \sum_{j=1}^{n+1} M_{i\leftarrow j}x_j(t),
\end{equation}
where $M_{i\leftarrow j}$ denotes the proportion of individuals of type $j$ that mutate to type $i$ in one generation. 
The existence of non-zero recombination between genes would require the addition of degree two terms in the polynomial map \cite{HS03},
while the effects due to genetic drift can be incorporated by means of 
diffusion approximations \cite{LH12}.

An understanding of the properties of maps from the simplex into itself is motivated by the study of these dynamical
systems. On the one hand, the deterministic dynamical systems mentioned above are  specified by one of these maps. On
the other hand, solving the diffusion partial differential equations (PDEs) that arise in the presence of genetic drift, 
requires the use of numerical methods that rely on the choice of a particular map from the simplex to itself\footnote{
For instance, in the case of the Finite Element Method, often one is interested in constructing an adaptive mesh on the 
simplex on which the PDE is defined. This is because the resolution needed to numerically solve the PDE is higher
in certain parts of the simplex than in others. A natural way to specify such an adaptive mesh consists of defining a uniform mesh
on the simplex, and an appropriate map from the simplex into itself. The choice of the map will be such that the preimages of the
subdomains where higher resolution is needed will span higher volumes than the preimages of the subdomains where lower resolution
is needed.}.

A particularly interesting class of maps are the polynomial maps of bounded-degree, which generalize the linear maps associated with
finite Markov chains to maps of arbitrary degree.
In the linear case, as represented in Eq. \eqref{blaone}, the general structure of the space of linear maps can be determined easily. 
In particular, as $\sum_{i=1}^{n+1}M_{i\leftarrow j}=1$ for each $j$, and every matrix element $M_{i\leftarrow j}$ is non-negative, 
the space of Markov matrices is the direct product of $n+1$ simplices of dimension $n$. 
We denote such a $n(n+1)$-polytope by $\mathsf{PMaps}(\Delta^{n},1)$, and by $\Delta^{n}$ we denote the $n$-dimensional unit simplex.
As the largest eigenvalue of $M$ is $1$ and because of the Perron-Frobenius theorem (see \cite{Sal97}), the dynamics of
any Markov chain in the interior of $\mathsf{PMaps}(\Delta^{n},1)$ always converges to a unique fixed point in the $n$-simplex.
This is because while one eigenvalue is $1$, the remaining eigenvalues are in the interior of the unit disk in the complex plane.
This compares with the Markov chains that are located in the boundary of the polytope $\mathsf{PMaps}(\Delta^{n},1)$,
where the particular arrangements of zeroes in the stochastic matrix span different faces.
In particular, stochastic matrices located in faces of decreasing 
dimensionality correspond to matrices with an increasing number of zeroes. 
The spectrum of these matrices consists of one eigenvalue equal to $1$, a set of eigenvalues
located in the interior of the unit disk in the complex plane, and another set of eigenvalues located in the boundary of the unit disk.
Indeed, when $M_{i\leftarrow j}$ has the maximum number of zero entries and its determinant is not zero, 
i.e. the matrix is a permutation matrix (see Appendix A for an elementary proof of this statement), all the eigenvalues
are located in the boundary of the unit disk. The spectrum of any of these stochastic matrices allows us to
determine some dynamical properties of the corresponding Markov chains; e.g. the existence of cycles is associated with more than one
eigenvalue in the unit circle of the complex plane (see for instance \cite{Sal97}).

In this paper, we are able to extend some of these known results in the linear case to the general case of 
the space of bounded-degree stochastic polynomial maps from the simplex into itself.
First, we demonstrate that these spaces are compact and convex subsets of Euclidean spaces. 
Second, we characterize a class of maps in the boundary of these spaces that generalize
the notion of stochastic matrix with maximum number of zeroes and nonzero determinant. 
Geometrically, these polynomial maps are defined by maximizing the number of preimages of the boundary of
the simplex. Indeed, this is satisfied by a family of maps that we introduce here, the \emph{minimal polynomial folding maps}, which, as we show below, are related to 
the theory of simplicial tessellations of the Euclidean space \cite{HW88, Ves91}.
Finally, we show how these maps can be constructed by solving certain systems of polynomial equations, we compute several examples
in one and two dimensions, and we study the dynamical properties of these maps and of their deformations in the space of polynomial maps.

\subsection{Basic notation and definitions.}
We denote by $\Delta^{n}\subset \mathbb{R}^{n+1}$ the standard $n$-simplex:
\begin{equation}
\Delta^{n}= \left\{ x=(x_1, x_2,\ldots, x_{n+1}) \in \mathbb{R}^{n+1}\colon\,\, \sum_{i=1}^{n+1} x_i=1, \,\, x_j\geq 0,\,  1\leq j \leq n+1  \right\}.
\end{equation}
For simplicity, the vertices, edges, faces and other boundary components of the $n$-simplex are denoted as its $0$-faces, 
$1$-faces, $2$-faces, \ldots, and $(n-1)$-faces. The interior of the $n$-simplex, $\Delta^{n}\backslash\partial \Delta^{n}$, is also referred to as
the $n$-face of the simplex. All the $i$-faces except for the $0$-faces, which are closed sets in $\mathbb{R}^{n+1}$, are open sets. Hence,
any simplex can be expressed as the disjoint union of its faces. This implies that if $\mathcal{S}$ is the set of all the faces in $\Delta^{n}$, there exists
an onto map $F_c\colon \Delta^n\to \mathcal{S}$, the \emph{face map}, which assigns to each point in the simplex its corresponding face, i.e., $x\in F_c(x)$.

Any polynomial map from the $n$-simplex to itself is defined by means of $n+1$ polynomials $\{ P_i(x_1,x_2,\ldots, x_{n+1}) \}_{i=1}^{n+1}$ 
satisfying
$$
\sum_{i=1}^{n+1} P_i(x)=1\quad \mathrm{and}\quad P_j(x)\geq 0,\, 1\leq j \leq n+1,
$$
for all $x\in \Delta^{n}$. We define the degree of a given map as the degree of the polynomial $P_i(x)$ of the largest degree. This definition
is ambiguous because any pair of polynomials, $P_i(x)$ and $P_i(x)(x_1+x_2+\cdots+x_{n+1})^N$ restricted to $\Delta^{n}$, are identical 
for any non-negative integer $N$. To avoid this ambiguity, we assume throughout that the degree of a polynomial $P(x)$ is the 
smallest degree possible after factoring out $(x_1+x_2+\cdots+x_{n+1})$ terms. 
We denote by $\mathsf{PMaps}(\Delta^{n},k)$ the space of polynomial maps of degree less than or equal to $k$ from $\Delta^{n}$ to itself. By
$\Pi_k^n$ we denote the vector space of polynomials in $n$ variables of degree at most $k$.

A $d$-folding map from the simplex to itself is a $d$-to-$1$ and onto continuous mapping in which the number of preimages is $d$
everywhere in the interior of $\Delta^{n}$. A polynomial fold of $\Delta^{n}$ is a particular instance
of a folding map that is also a polynomial map. 

\section{The geometry of the space of stochastic polynomial maps}

For a given positive integer $d$, the corresponding $d$-folding polynomial maps located in the space of polynomial maps
$\mathsf{PMaps}(\Delta^{n},k)$ with the lowest degree $k$ possible, are a very special set of points
in $\mathsf{PMaps}(\Delta^{n},k)$. 
This will become clearer after we describe in more detail the geometry of $\mathsf{PMaps}(\Delta^{n},k)$.

\begin{theorem} 
The space of polynomial maps of degree $k$ from the $n$-simplex to itself, $\mathsf{PMaps}(\Delta^{n},k)$, is a compact
and convex subset of the vector space of polynomial maps of degree $k$ from $\mathbb{R}^n$ to itself.
\end{theorem}
\emph{Proof}. The proof consists of an explicit construction of $\mathsf{PMaps}(\Delta^{n},k)$, such that the
properties stated in the theorem follow from such construction.
For simplicity, we work with the projection of the simplex $\Delta^n\subset \mathbb{R}^{n+1}$ down to $\mathbb{R}^{n}$, where
$\Delta^{n}$ is defined by the set of inequalities
$$
x_i\geq 0,\,\, 1\leq i\leq n,
$$
$$
\sum_{i=1}^n x_i\leq 1.
$$
Any polynomial map $f$ in $\mathsf{PMaps}(\Delta^{n},k)$ can be defined by means of $n$ polynomials $\{ P_i(x) \}_{i=1}^n$ of degree
less than or equal to $k$, that obey the following conditions:
\begin{eqnarray}
P_i(x) \geq 0, \,\, \forall i\in [1,n],\, \forall x\in\Delta^n, \label{firstcondition} \\
1-\sum_{i=1}^n P_i(x) \geq 0,\,\,\forall x\in\Delta^n,\label{secondcondition}
\end{eqnarray}
and 
\begin{equation}
x^\prime=f(x), \quad x^\prime_i=P_i(x), \label{thirdcondition}
\end{equation}
which states how preimages $x$ are mapped to images $x^\prime$ and therefore, it completes the definition of $f$. 
The condition in Eq. \eqref{firstcondition} can be simplified thanks to an important theorem of P\'olya, which provides a systematic process for
deciding whether a given polynomial is strictly positive on the simplex \cite{HLP52}. Using P\'olya's result, one can show how
the polynomials that are non-negative on $\Delta^n$
span a convex cone in $\Pi_k^n$, the space of polynomials in $n$ variables of degree at most $k$. To show this, we first need
to construct homogeneous representatives of the polynomials in $\Pi_k^n$.
In particular, if $P(x)$ is 
a polynomial in $\Pi_k^n$, we construct its homogenous polynomial representative 
$P_H(x)$ in $\mathbb{R}^{n+1}$ as follows. We write $P(x)$ as a sum of monomials
$$
P(x)=\sum_{\vert \alpha \vert\leq k} c_{\alpha_1 \alpha_2\cdots\alpha_n}x^{\alpha_1}_1x^{\alpha_2}_2\cdots x^{\alpha_n}_n
$$
with $\vert \alpha \vert=\sum_{i=1}^n \alpha_i$, and multiply the terms with $\vert \alpha \vert < k$ by $(x_1+x_2+\cdots +x_{n+1})^{k-\vert \alpha \vert}$.
In this way, we obtain a homogenous polynomial of degree $k$ which is identical to $P(x)$ when restricted to $\Delta^n\subset\mathbb{R}^{n+1}$: $P_H(x)=$
$$
\sum_{\vert \alpha \vert\leq k} c_{\alpha_1 \alpha_2\cdots\alpha_n}\left( \delta_{\vert \alpha \vert,k}+ (1-\delta_{\vert \alpha \vert,k})(x_1+x_2\cdots +x_{n+1})^{k-\vert \alpha \vert} \right) x^{\alpha_1}_1x^{\alpha_2}_2\cdots x^{\alpha_n}_n,
$$
with $\delta_{\vert \alpha \vert,k}=1$ when $\vert \alpha \vert=k$ and $\delta_{\vert \alpha \vert,k}=0$ otherwise. Using
this homogenous polynomial $P_H(x)$, P\'olya's theorem simply says that if $P(x)$ is positive on $\Delta^n$, i.e. $P(x)>0$, $\forall x\in\Delta^n$,
then there exists a non-negative integer $N$ such that all the coefficients of $(x_1+x_2+\cdots +x_{n+1})^N P_H(x)$ are positive.
If $P(x)$ is positive on $\Delta^n$ then $\lambda P(x)$, with $\lambda\in\mathbb{R}^+$, is also positive. This implies that the 
polynomials of degree $k$ with $P(x)>0$, $\forall x\in\Delta^n$, span an open convex cone $\mathcal{K}_{+}(\Delta^n,k)\subset \Pi_k^n$.  
It also follows that $\dim\,\mathcal{K}_{+}(\Delta^n,k)=\dim\, \Pi_k^n$. Furthermore, one can use this construction
to describe the set of non-negative polynomials on the simplex. More precisely,
if $P(x)$ is strictly positive on the $n$-simplex and $Q(x)$ is non-negative for all $x\in\Delta^n$, then any deformation of $Q(x)$ by $\epsilon P(x)$
makes $Q(x)+\epsilon P(x)$ an element of $\mathcal{K}_{+}(\Delta^n,k)$ for any positive real value of $\epsilon$.
Therefore, the closure of $\mathcal{K}_{+}(\Delta^n,k)$ in $\Pi_k^n$, which we denote by $\mathcal{K}(\Delta^n,k)=\mathrm{cl}\left( \mathcal{K}_{+}(\Delta^n,k) \right)$, is the set of non-negative polynomials on $\Delta^n$. In other words,
$$
\{ P(x) \in \Pi_k^n\colon\,\, P(x)\geq 0,\,\, \forall x\in \Delta^n \}=\mathcal{K}(\Delta^n,k).
$$
Thus, the set of polynomials $\{ P_i(x)\}_{i=1}^n$ that satisfy Eq. \eqref{firstcondition} are elements of the cartesian product
of $n$ closed convex cones $\mathcal{K}(\Delta^n,k)$. This product of cones is itself a closed convex cone in $\oplus_{i=i}^n \Pi_k^n$,
the direct sum of vector spaces of polynomials on the simplex. We denote such a convex cone by $\times_{i=1}^n \mathcal{K}(\Delta^n,k)$.

The space $\mathsf{PMaps}(\Delta^{n},k)\subset \oplus_{i=i}^n \Pi_k^n$ and its compactness property arise after we impose Eq. \eqref{secondcondition} 
on the convex cone $\times_{i=1}^n \mathcal{K}(\Delta^n,k)$. By the Heine-Borel theorem it is sufficient
to show that $\mathsf{PMaps}(\Delta^{n},k)$ is a bounded subset of the cone $\times_{i=1}^n \mathcal{K}(\Delta^n,k)\subset \oplus_{i=i}^n \Pi_k^n$, 
in order to demonstrate that $\mathsf{PMaps}(\Delta^{n},k)$ is compact. First, let 
$$\mathrm{span}_\lambda \{ \lambda P_i(x) \}_{i=1}^n$$ 
for $\lambda\in\mathbb{R}^+$ positive and real, denote the ray in $\times_{i=1}^n \mathcal{K}(\Delta^n,k)$ generated by $\{ P_i(x) \}_{i=1}^n$. The maximum
of the function $\sum_{i=1}^n P_i(x)$ on the simplex is strictly positive for any non-zero set of polynomials $\{ P_i(x) \}_{i=1}^n \in 
\times_{i=1}^n \mathcal{K}(\Delta^n,k)$. 
Let us denote such a maximum value by $S_{max}$. Now, Eq. \eqref{secondcondition} restricted on this ray is equivalent to requiring the non-negativity of
$1-\lambda \sum_{i=1}^n P_i(x)$ for all $x$ in the simplex. This is only satisfied by 
the segment $\{ \lambda P_i(x) \}_{i=1}^n$ parametrized by $\lambda\in [0,1/S_{max}]$. Similarly, for every ray in
$\times_{i=1}^n \mathcal{K}(\Delta^n,k)$ only a bounded and finite subinterval of the ray satisfies Eq. \eqref{secondcondition}.
It follows that $\mathsf{PMaps}(\Delta^{n},k)$ is a bounded and closed subset of the convex cone $\times_{i=1}^n \mathcal{K}(\Delta^n,k)$. 

The convexity of $\mathsf{PMaps}(\Delta^{n},k)$ follows from the convexity of the simplex $\Delta^n$. In particular, if
the image of $x$ defined by $\{ P_i(x) \}_{i=1}^n \in \mathsf{PMaps}(\Delta^{n},k)$ is in $\Delta^n$ for each $x\in\Delta^n$, 
and the same is true for a second map $\{ Q_i(x) \}_{i=1}^n\in \mathsf{PMaps}(\Delta^{n},k)$, then $\{ tP_i(x) +(1-t)Q_i(x) \}_{i=1}^n$
is in $\Delta^n$ and yields a well-defined map in $\mathsf{PMaps}(\Delta^{n},k)$ for any $t\in [0,1]$.$\square$
\vspace{0.75cm}

\begin{figure}[hpbt]
\begin{center}
\vspace{0.5cm}
\begin{tabular}{c}
    \includegraphics[width=10.0cm]{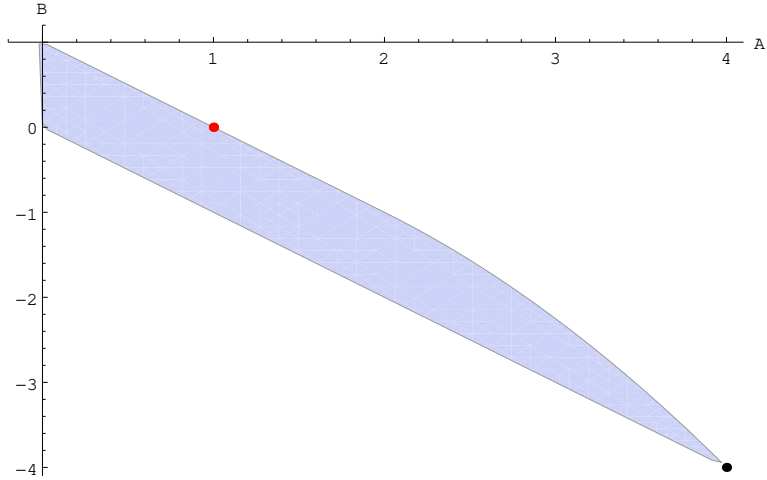}\\
\end{tabular}
\caption{\label{figc} Convex set of quadratic maps of the type $f(x)=x(A+Bx)$ on the $1$-simplex.
The horizontal axis parametrizes the $A$-coefficient and the vertical axis parametrizes the $B$-coefficient. For
pairs $\{A,\, B\}$ outside the convex set, the image of the quadratic map is not in the simplex.
The red dot denotes the identity map, while the black dot denotes the polynomial two-folding map
$f_2(x)=\frac{1-\cos(2\arccos (1-2x))}{2}=4x(1-x).$}
\end{center}
\end{figure}

Our construction of $\mathsf{PMaps}(\Delta^{n},k)$ relies on P\'olya's criterion to decide whether a given polynomial 
vanishes somewhere in the simplex $\Delta^n$ or not. In particular, we have showed how the interior of $\mathsf{PMaps}(\Delta^{n},k)$ 
consists of maps for which every defining polynomial, $\{ P_i(x)\}_{i=1}^n$ plus the polynomial $1-\sum_j P_j(x)$, are strictly positive on $\Delta^n$, and therefore, 
their corresponding images are in the interior of $\Delta^n$. This compares with the boundary of 
$\mathsf{PMaps}(\Delta^{n},k)$, consisting of maps for which some of the defining polynomials 
do vanish on some algebraic subsets of the standard simplex. For these maps the preimages of the boundary of $\Delta^{n}$ 
are not the empty set.
The case of finite Markov chains is particularly illuminating in this regard.
When $k=1$, the defining polynomials are non-negative linear forms on the simplex, while the boundary
of $\mathsf{PMaps}(\Delta^{n},1)$ consists of maps that have a non-empty preimage of the boundary $\partial \Delta^n$.
In an extreme limiting case, there exist maps with non-singular Jacobian matrices that maximize the size 
of the corresponding preimages of $\partial \Delta^n$. These maps consist of invertible stochastic matrices that have the largest number possible of zero entries, i.e. the permutation
matrices. Geometrically, these matrices are special vertices in the boundary of the
$n(n+1)$-dimensional polytope $\mathsf{PMaps}(\Delta^{n},1)$, among other things because they are the only
bijective maps therein (see Lemma 1 in the Appendix).

One can define nonlinear analogs of these extremal maps in $\mathsf{PMaps}(\Delta^{n},k)$ for $k>1$. They consist of 
polynomial maps $f\colon\, \Delta^n\to\Delta^n$ with a non-singular Jacobian almost everywhere in $\Delta^n$, that
maximize the complexity of the algebraic subset in $\Delta^n$ that is the preimage of the boundary under $f$, i.e. they maximize the complexity
of $f^{-1}\left( \partial \Delta^n \right)\subset \Delta^n$. These maps are easy to construct in the case of the $1$-simplex.  
For instance, when $n=1$ and $k=2$ the largest algebraic subset of the interval that can be mapped to its boundary by a quadratic map
(with a non-singular Jacobian almost everywhere in $[0,1]$) consists of three points (see Fig. 1). In general,
the largest algebraic subset of $\Delta^1$ that can be mapped to the boundary of $\Delta^1$ by a map in $\mathsf{PMaps}(\Delta^{1},k)$ consists
of $k+1$ points. These maps can be constructed explicitly \cite{HW88,Ves91}. They coincide with the Chebyshev maps
and their permutations:
\begin{equation}
f_k\colon\, [0,1]\to [0,1],\quad f_k(x)=\frac{1-T_k(1-2x)}{2}.
\end{equation}
Here, $T_k$ denotes the $k^{th}$ Chebyshev polynomial of the first kind, which can be evaluated as $T_k(x)=\cos\left(k\arccos x \right)$.
One can show from the properties of the Chebyshev maps that every map $f_k$ is a $k$-folding map of the $1$-simplex.

As far as we know, there do not exist previous attempts in the literature to construct 
the higher dimensional analogs of these extremal maps on the unit simplex.
Although these maps are indeed harder to construct, we show below how some particular families of 
extremal maps admit a simple geometrical description. 
In particular, the concept of a polynomial
map $f$ of fixed degree that maximizes the complexity of $f^{-1}\left( \partial \Delta^n \right)$ finds an instance
in any polynomial $d$-folding map of minimal degree. 
As every polynomial $d$-folding map $f$ is continuous, onto and $d$-to-$1$ everywhere in the interior of the image simplex, one can decompose the preimage simplex 
$\Delta^n$ as the union of $d$ sets, $\Delta^n=\cup_{a=1}^d \delta_a$, such that the restriction of the map on each subset,
$f\vert_{\delta_a}\colon\, \delta_a\to\Delta^n$, is one-to-one and onto. In this case, the preimage of the boundary is simply the union
of the boundaries of each subset $\delta_a$,
$f^{-1}\left( \partial \Delta^n \right)=\cup_{a=1}^d \partial \delta_a$. Therefore, those $d$-folding maps $f$ located in spaces
$\mathsf{PMaps}(\Delta^{n},k)$ that have the smallest $k$ possible maximize the complexity of $f^{-1}\left( \partial \Delta^n \right)\subset \Delta^n$
among all maps in $\mathsf{PMaps}(\Delta^{n},k)$.
In the next section we define and construct these minimal polynomial folding maps. First, we propose a restricted notion of homotopy
that allows us to group together equivalent folding maps. For instance, we consider $4x(1-x)$ and $2x\Theta(1/2-x)+(2-2x)\Theta(x-1/2)$
to be topologically equivalent $2$-folding maps of the $1$-simplex (here, $\Theta$ denotes the Heaviside function); however, 
we do not consider the polynomial $2$-fold $4x(1-x)$ and the piecewise linear $3$-fold $3 x \Theta(1/3-x)+(2-3 x)\Theta(x-1/3)\Theta(2/3-x)+(3 x-2)\Theta(x-2/3)$ 
to be equivalent. Then, we use this notion of homotopy
to define the \emph{minimal polynomial folding maps}. Finally, we show how to construct these maps as the zero loci 
of certain polynomial systems equations and compute several examples.

\section{Polynomial Folding Maps}

In order to be able to compare polynomial folds with non-polynomial folds, e.g. piecewise linear folding maps, we need an
appropriate criterion to establish a topological equivalence relation between them.
To that end, we introduce a restricted notion of homotopy equivalence between folding maps as follows.
\begin{figure}[hpbt]
\begin{center}
\vspace{0.5cm}
\begin{tabular}{cc}
    \includegraphics[width=6.0cm]{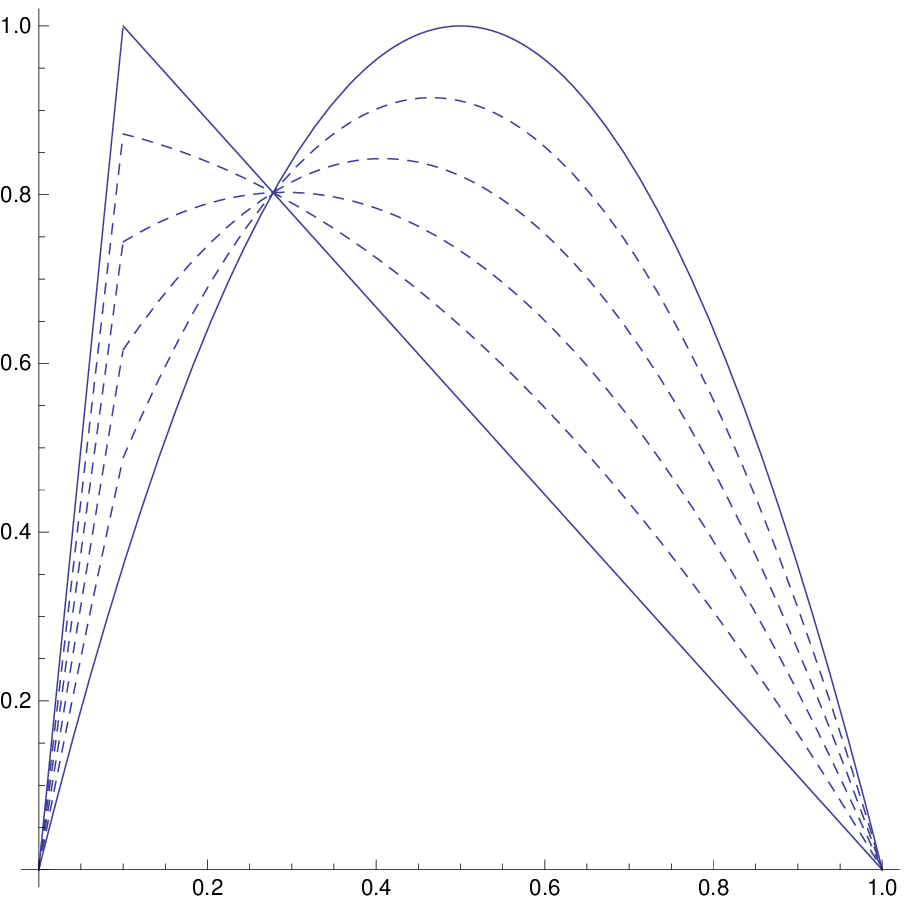}&
    \includegraphics[width=6.0cm]{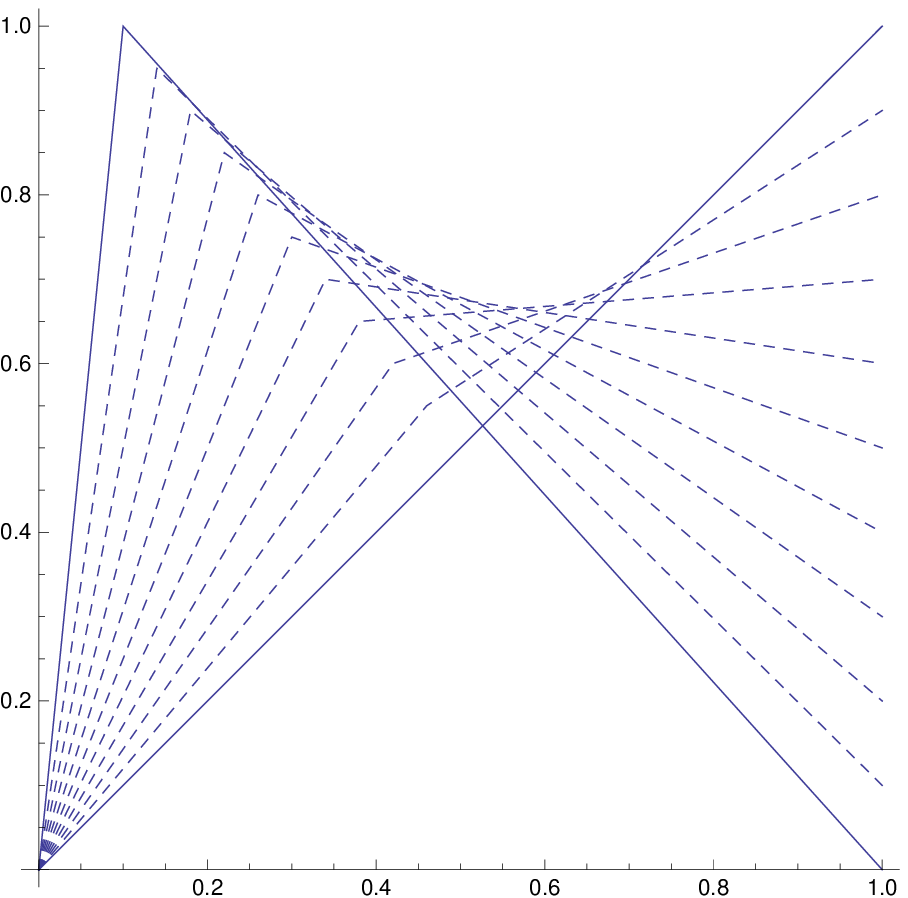}\\
\end{tabular}
\caption{\label{figa} A piecewise linear $2$-fold of the $1$-simplex. Here, the $1$-simplex 
spans the interval $[0,1]$, the horizontal axis denotes the preimage of the map, and the vertical axis denotes its
image. The standard notion of homotopy equivalence between maps is not valid to compare folding maps.
The $2$-folding map, given by the function $f(x)=10 x\Theta(\frac{1}{10}-x)+(1-
\frac{10}{9}(x-\frac{1}{10}))\Theta(x-\frac{1}{10})$, is homotopic to the polynomial fold $g(x)=4x(1-x)$ (see left plot). The same $2$-folding
map is also homotopic to the identity map $I(x)=x$, (see right plot).
}

\end{center}
\end{figure}
We say that two maps, $f\colon\, \Delta^{n}\to\Delta^{n}$ and $g\colon\, \Delta^{n}\to\Delta^{n}$,
are \emph{face-constrained homotopic} (or \emph{fc-homotopic}) if and only if there exist two continuous functions, 
$H\colon\, \Delta^{n}\times[0,1]\to \Delta^{n}$ and $K\colon\, \Delta^{n}\times[0,1]\to \Delta^{n}$,
from the product of the $n$-simplex with the unit interval $[0,1]$ to the $n$-simplex such that:
\begin{itemize}
\item[i)] $K(x,t)\colon\, \Delta^{n}\to\Delta^{n}$ is one-to-one and onto for all $t\in [0,1]$.
\item[ii)] $\{ K(x,0),\, H(x,0)\}=\{ x,\, f(x)\}$ and $\{ K(x,1),\, H(x,1)\} =\{ y,\, g(y)\}$ for all $x \in \Delta^{n}$.
\item[iii)] The \emph{face-constrained} condition requires that if $F_c(x)$ is the face
where $x$ belongs, $x\in F_c(x)\subset \Delta^n$, and $F_c\left( f(x)\right)$ denotes the face where the image $f(x)$ is located, $f(x)\in F_c\left( f(x)\right)\subset \Delta^n$, then the deformation defined by $K$ and $H$ preserves such constraints, i.e. $K(x,t)\in F_c(x)$ and $ H(x,t)\in F_c\left( f(x)\right)$ for all $t\in [0,1]$.
\end{itemize}
Note that if the condition iii) is not satisfied, one recovers the standard notion of homotopy by choosing
$K$ to be the identity map.

\begin{figure}[hpbt]
\begin{center}
\vspace{0.5cm}
\begin{tabular}{c}
    \includegraphics[width=6.0cm]{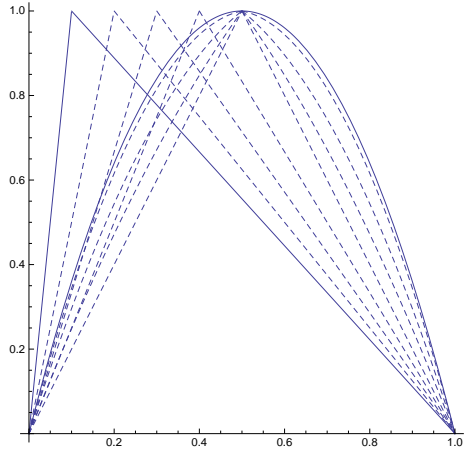}\\
\end{tabular}
\caption{\label{figb} A \emph{face-constrained homotopy} of the piecewise linear fold defined by $f(x)=10 x\Theta(\frac{1}{10}-x)+(1-
\frac{10}{9}(x-\frac{1}{10}))\Theta(x-\frac{1}{10})$. The preimage points $x=0$ and $x=1$ are located on $0$-faces, while both of their images 
are on the same $0$-face ($y=0$). The images of all the interior points, except $x=1/10$, are also interior points. The image of the interior point $x=1/10$
is located on the $0$-face $y=1$. These constraints are preserved by any \emph{fc-homotopy}. Hence, while $f(x)$ is \emph{fc-homotopic} to $g(x)=4x(1-x)$, it is not \emph{fc-homotopic} to the identity map.}
\end{center}
\end{figure}

\begin{definition}
\emph{Minimal polynomial folding map}. Given an fc-homotopy class of folding maps on the $n$-simplex, we say that 
a particular polynomial representative is minimal if and only if the corresponding polynomial map exists in $\mathsf{PMaps}(\Delta^{n},k)$ 
and there do not exist fc-homotopically equivalent maps in $\mathsf{PMaps}(\Delta^{n},l)$, with $l<k$.
\end{definition}
Let us assume that a given minimal polynomial folding map $f\colon\, \Delta^n\to\Delta^n$ is defined by the polynomials
$\{P_i(x)\}_{i=1}^n$ and $P_{n+1}(x)=1-\sum_{j=1}^n P_j(x)$. 
Then, the preimage of each of the $n+1$ facets of the simplex is an algebraic subset of the unit simplex defined 
either by 
$$
\{ x\in \Delta^n \colon\, P_i(x)=0 \}, 
$$ 
for the facet $y_i=0$, or by
$$
\{ x\in \Delta^n\colon\, 1-\sum_{j=1}^n P_j(x)=P_{n+1}(x)=0 \},
$$
for the facet $\sum_{i=1}^n y_i=1$.
We denote by $\{ z_i \subset \Delta^n\backslash \partial\Delta^n ,\, b_i \subset \partial\Delta^n \}$ the codimension-one preimage of the $i^{th}$ facet, decomposed in its interior and boundary parts.
In this case, $z_i$ is the subset of the preimage located in the interior of the simplex and 
$b_i$ is the complementary subset located in the boundary of the simplex. 
It is important to note that as
$b_i\subset \partial\Delta^n$ is the zero locus of a polynomial defined on $\Delta^n$, $b_i$ is either the empty set, or a finite set of facets.
Using other words, $b_i$ cannot be composed of finite subsets within facets because such subsets are not algebraic; 
$b_i$ can only be composed of whole facets, because they are the only algebraic subsets of codimension one in the boundary of the simplex, 
e.g. $x_i=0$, or $1-\sum_{j=1}^n x_j=0$. Indeed, there exist folding maps
that do not admit fc-homotopic polynomial representatives; this is because whenever
different subsets within a single facet are mapped to distinct facets, the folding map cannot be polynomial.
This will become clearer after we show the following.

\begin{theorem}
Let $f\in \mathsf{PMaps}(\Delta^{n},k)$ be a minimal polynomial folding map of degree $k$.
Then, each defining polynomial $P_i(x)$, with $1\leq i\leq n+1$, can be factored as
\begin{equation}
\label{eq_factor}
P_i(x)=l_i(x)(q_i(x))^2 m_i(x),
\end{equation}
where $l_i(x)$, $q_i(x)$ and $m_i(x)$ obey the following:
\begin{itemize}
\item The zero locus of $l_i(x)$ in $\Delta^n$ is exactly $b_i\subset \partial\Delta^n$. 
\item The degree of $l_i(x)$ equals the number of facets contained in $b_i$.
Furthermore, the product of all the $l_i(x)$'s is
$$
\prod_{i=1}^{n+1}l_i(x)=(1-\sum_j x_j)\prod_{i=1}^n x_i.
$$
\item The zero locus of $q_i(x)$ in the interior of the $n$-simplex is exactly $z_i \subset \Delta^n\backslash \partial\Delta^n$.
\item $l_i(x)$ and $m_i(x)$ are strictly positive in the interior of the $n$-simplex. The zero set of $m_i(x)$ is either empty or 
it is contained in the zero set of $l_i(x)$.
\item The degree of $P_i(x)$ is bounded by $k$, $\deg l_i +2\deg q_i +\deg m_i \leq k$.
\end{itemize}
\end{theorem}
\emph{Proof}.
The preimage of $\Delta^n\backslash\partial\Delta^n$ under $f$ consists of
$d$ disconnected open subsets in the interior of $\Delta^n$. The restriction
of $f$ to any of these subsets gives rise to a one-to-one mapping that has the whole interior of the simplex as codomain. 
The complementary set of the union of these $d$ open subsets in $\Delta^n$ is $f^{-1}(\partial\Delta^n)$.
These facts follow from the properties of the folding maps:
they are continuous and onto, for every point in the interior of the simplex there exist $d$ different
preimages, and they map the boundary $\partial\Delta^n$
to a subset in the boundary of the image simplex (see Lemma 2 in the Appendix). In addition to this, 
the fact that a minimal polynomial folding map is also polynomial imposes further constraints. In particular,
the properties of the defining polynomials stated in Theorem $2$ are derived from the algebraic properties 
of the preimage set $f^{-1}(\partial\Delta^n)$. 
Algebraically, $f^{-1}\left( \partial\Delta^n \right)$ is the union of irreducible varieties  
of codimension one in $\Delta^n$. The set defined by the zero locus $P_{i}(x)=0$ denotes the 
algebraic variety in $f^{-1}\left( \partial\Delta^n \right)$ that is preimage of the $i^{th}$ facet.
As such a preimage has to be the union of a finite number of facets and a codimension one
variety in the interior of the simplex, $P_i(x)$ can be factored as
$$
P_i(x)=l_i(x)r_i(x)m_i(x).
$$
Here, the zero locus of $l_i(x)$ is the union of facets in the boundary,
the zero locus of $r_i(x)$ is a codimension one variety in the interior of the simplex, and $m_i(x)$
is a non-negative polynomial of degree equal or less than $k-\deg l_i(x)- \deg r_i(x)$, whose zero locus is either empty or it is contained 
in the zero locus of $l_i(x)$.
The existence of $m_i(x)\in\mathcal{K}(\Delta^n,$
$k-\deg l_i(x)- \deg r_i(x))$ follows from the facts that $P_i(x)$ has a degree equal or less than $k$, 
that $P_i(x)\geq 0$ for all $x\in\Delta^n$ and that the zero set of $P_i(x)$ is already 
specified by $r_i(x)$ and $l_i(x)$. Because of these constraints, the existence of an additional polynomial factor of
degree less than or equal to $k-\deg l_i(x)- \deg r_i(x)$ is possible. 

As the $i^{th}$-facet of the $n$-simplex is defined by the zero
locus of the degree-one monomial $x_i$ (when $i=n+1$, $x_{n+1}=1-\sum_{j=1}^n x_j$), the lowest degree polynomial 
whose zero locus consists of a predefined set of facets $b_i\subset \Delta^n$ is
$$
l_i(x)=\prod_{j\in b_i} x_j.
$$
If $b_i$ is empty then $l_i(x)=1$. It follows that $l_i(x)>0$ for all $x$ in the interior of the simplex.
Furthermore, as $f(\partial\Delta^n)$ is contained in $\partial\Delta^n$, the images of every facet in $\Delta^n$
are in the union of all the image facets, and therefore
$$
\prod_{i=1}^{n+1} l_i(x)
$$
has to be a polynomial divisible by $(1-\sum_{j=1}^nx_j)\prod_{i=1}^n x_i$. As the preimage of each facet is the union
of irreducible varieties of codimension one, each facet has to be mapped to a whole facet. This implies that the monomial $x_j$
whose zero locus defines the facet cannot appear as a factor in more than one polynomial $l_i(x)$. Otherwise, the image of such  $j^{th}$-facet
would be a face of codimension larger than one. It follows that $\prod_{i=1}^{n+1} l_i(x)$
is not only divisible by $(1-\sum_{j=1}^nx_j)\prod_{i=1}^n x_i$ but it is a constant multiple of $(1-\sum_{j=1}^nx_j)\prod_{i=1}^n x_i$.
By convention, such a constant factor is absorbed in the $m_i(x)$ terms, and 
$$
\prod_{i=1}^{n+1} l_i(x)=(1-\sum_{j=1}^nx_j)\prod_{l=1}^n x_l.
$$

Finally, it remains to show that the factor $r_i(x)$ is the square of a polynomial: $r_i(x)=\left( q_i(x) \right)^2$.
As the zero set $z_i=$ $\{x\in\Delta^n\colon\, r_i(x)=0\}$ is a codimension-one variety in the interior of the simplex, 
and $r_i(x)$ is strictly positive for every $x\in\Delta^n\backslash z_i$, it follows from the Taylor expansion of $r_i(x)$ around any $x\in z_i$ that
the gradient of $r_i(x)$ restricted to $z_i$ has to vanish. Therefore, $z_i$ has to be the 
zero locus of a polynomial $q_i(x)$ which might have negative values on a subset of the simplex and satisfies $r_i(x)=(q_i(x))^2$.$\square$
\vspace{0.75cm}

A priori, the factorization in Eq. \eqref{eq_factor} can be used 
to construct any minimal polynomial folding map. First,
one fixes an fc-homotopy class of folding maps that admits polynomial representatives. Then, the assignment
of facets to image facets defined by the fc-homotopy class is used to determine the polynomials $\{ l_i(x) \}_{i=1}^{n+1}$.
Unfortunately, specifying
the remaining factors in Eq. \eqref{eq_factor} is a harder problem. The approach that we follow consists of
parametrizing the set of all the polynomials, $\{ m_i(x,\phi^i) \}_{i=1}^{n+1}$ and $\{ q_i(x, \theta) \}_{i=1}^{n+1}$, that are
compatible with the given class of folding maps. Then, we determine the particular parameters $\{\hat{\theta},\, \hat{\phi}^i \}_{i=1}^{n+1}$ 
that define the folding map by solving a system of polynomial equations. More precisely, we use the defining equation
\begin{equation}
\label{defeq}
l_{n+1}(x)\left( q_{n+1}(x, \theta) \right)^2 m_{n+1}(x,\phi^{n+1}) = 1-\sum_{i=1}^n l_{i}(x)\left( q_{i}(x,\theta) \right)^2 m_{i}(x,\phi^i),
\end{equation}
to determine the parameters $\{\hat{\theta},\, \hat{\phi}^i \}_{i=1}^{n+1}$ that satisfy Eq. \eqref{defeq} for all $x\in\Delta^n$. 
In practice, this is a tedious approach 
because the set of compatible polynomials usually has several independent components and each component is a complicated semialgebraic
set in a high-dimensional vector space. Nevertheless, it is possible to identify each component and construct the corresponding spaces
of polynomials.
In particular, every component is associated with a different partition $\{(a_i,b_i)\}_i$ of the degrees of $(q_i(x))^2m_i(x)$; 
i.e. partitions of $k-\deg l_i(x)=2a_i+b_i$,
with $a_i=\deg q_i(x)$, $b_i=\deg m_i(x)$ and $k$ the degree of the polynomial map. 
Given a particular partition $\{(a_1,b_1),$ $(a_2,b_2),$ $\ldots (a_{n+1},b_{n+1}) \}$
of the degree of every polynomial factor, it is possible to construct the space of compatible
polynomial factors with these prescribed degrees and find the solutions of Eq. \eqref{defeq}.

This construction becomes more transparent in the next subsection, where we show several examples in one and two dimensions.
However, before we construct particular examples, it is useful to state some general properties of these families of polynomials.
For instance, the family of polynomials $\{ q_i(x, \theta) \}_{i=1}^{n+1}$ can be determined by requiring
the corresponding zero loci  $\{ z_i\subset \Delta^n \}_{i=1}^{n+1}$ to be compatible with the given class of folding maps. 
This allows for the possibility of the existence of maps in the fc-homotopy class
for which the preimage of $\partial \Delta^n$ restricted to the interior of the simplex is exactly defined by the set of 
hypersurfaces $\{ z_i \}_{i=1}^{n+1}$. The parameter space $\{ \phi^i \}$ admits a simpler description.
The set of compatible polynomials $m_i(x, \phi^i)$  consists of those polynomials of degree $b_i$ 
$$
m_i(x, \phi^i)= \sum_{\vert \alpha\vert\leq b_i} \phi^{i}_{\alpha_1 \alpha_2\cdots\alpha_n}x_1^{\alpha_1} x_2^{\alpha_2} \cdots x_n^{\alpha_n}
$$
whose zero loci are either non-existent or are contained in the zero locus of $l_i(x)$.

\subsection{Examples of minimal polynomial folding maps}

In the following examples we fix the fc-homotopy class where the minimal polynomial folding map sits by defining a 
non-polynomial folding map in the same class of maps.
In particular, a well studied family of folding maps are the piecewise linear folding maps derived from periodic tessellations of the $n$-dimensional
Euclidean space \cite{HW88}. 
In this construction one defines an embedding of the simplex into a Euclidean space of the same
dimension, such that the action of the group of reflections through the facets of the simplex generates a tessellation of the space 
(e.g. see Fig. 4). One then uses this tessellation to determine the associated folding maps. More precisely,
we can construct these piecewise linear folding maps by means of (see \cite{HW88} for more details on this construction):
\begin{itemize}
\item[i)] An affine embedding of $\Delta^n$ in $\mathbb{R}^n$ such that the embedded simplex tessellates $\mathbb{R}^n$ by means of
reflections through the hyperplanes associated with the $(n-1)$-faces of the simplex. 
We denote this affine map by $A\colon\, \Delta^n\to\mathbb{R}^n$, and we 
assume that the image of one of the vertices of $\Delta^n$ is the zero vector in $\mathbb{R}^n$. 
\item[ii)] An integer dilation of the image simplex in $\mathbb{R}^n$. In particular, we require that for any integer dilation, 
the dilated simplex is tessellated by a finite number of smaller replicas of itself (e.g. see Fig. 4.).
\item[iii)] A map that assigns to each point in every tile of the tessellation of $\mathbb{R}^n$ the equivalent point in the 
tile that is the image of $\Delta^n$ under $A$. 
This map is unique; it is constructed by means of reflections through hyperplanes; and we denote it by
$h\colon\, \mathbb{R}^n\to A(\Delta^n)$.
\end{itemize}

Now, we define a folding map by first dilating the affine simplex $A(\Delta^n)$ $m$ units, 
i.e. we multiply every vector in $A(\Delta^n)\subset\mathbb{R}^n$ by the positive integer $m$.
We denote such dilation by $m\circ A(\Delta^n)$, and by
$d$ we denote the number of smaller replicas of the simplex contained in $m\circ A(\Delta^n)$. Therefore, the piecewise linear $d$-folding map is
defined as the following composition of mappings 
$$
A^{-1}\circ h \circ m\circ A\colon\, \Delta^n\to\Delta^n.
$$
In what follows we define several piecewise linear folding maps in one and two dimensions, and then apply the polynomial factorization 
of Theorem 2 to determine their associated minimal polynomial folding maps.

\subsubsection{Folding maps of the one-simplex}

The piecewise linear $d$-folding maps in dimension one can be derived from partitions of the unit interval in $d$ subsegments. 
In particular, if we assume that $x=0$ is mapped to $0$ then a partition determines uniquely a 
piecewise linear folding map. The simplest non-trivial case, the two-folding map, can be written as 
$f_{2,a}(x)=(x/a)\Theta(a-x)+((1-x)/(1-a))\Theta(x-a)$ with $a\in (0,1)$; although for simplicity we consider only $a=1/2$. 
It follows from the definition of $f_{2,a}$ that there does not exist a strictly linear map that is fc-homotopic to this two-folding map. 
However, if there exists a polynomial map of degree two that is fc-homotopic to the two-fold, then by Theorem 2 its corresponding 
defining polynomials will factor as
\begin{equation}
\label{dimonefold}
P_1(x)=Ax(1-x), \quad P_2(x)=1-P_1(x)=(B+Cx)^2,
\end{equation}
with $A>0$, $B$ and $C$ being unknowns to be determined. The fc-homotopy equivalence to $f_{2,a}(x)$ requires that the preimages of $y=0$ are $x=0$ and $x=1$,
and that there exists an interior point $x_i=-B/C\in (0,1)$ that is mapped to $y=1$. We can determine
the coefficients $A$, $B$ and $C$ by expanding the defining equation $P_2(x)=1-P_1(x)$ in powers of $x$
$$
(B+Cx)^2=1-Ax(1-x), \quad B^2-1+(2BC+A)x+(C^2-A)x^2=0.
$$
If this is true for every $x\in [0,1]$, then $B=1$, $C=-A/2$, $A^2/4=A$, and $A=4$. Therefore, the corresponding minimal
polynomial folding map is defined by $P_1(x)=4x(1-x)$ and $P_2(x)=(1-2x)^2$, i.e. the logistic map \cite{May76}.

\begin{table}[ht]
\begin{center}
\begin{tabular}{|c c c|}
\hline\hline
$k$ & $P_1(x)$ & $P_2(x)=1-P_1(x)$\\ [0.5ex]
\hline\hline
$1$ & $x$ & $1-x$ \\
\hline
$2$ & $4x(1 - x)$ & $(1-2x)^2$ \\ 
\hline
$3$ & $x (3 - 4 x)^2$ & $(1-x)(1-4 x)^2$ \\
\hline
$4$ & $16x(1-x) (1 - 2 x)^2$ & $\left(8 x^2-8 x+1\right)^2$ \\
\hline
$5$ & $x (5 - 20 x + 16 x^2)^2$ & $(1-x) \left(16 x^2-12 x+1\right)^2 $ \\
\hline
$6$ & $4x(1-x)\left( (1 - 4 x)(3 - 4 x) \right)^2$ & $\left( (1 - 2 x) (1 - 16 x + 16 x^2)\right)^2$ \\
[1ex]
\hline
\end{tabular}
\caption{\label{tableA} Polynomial factorization of the first six Chebyshev maps from the 
$1$-simplex to itself.}
\end{center}
\end{table}

One can repeat this construction for every fc-homotopy class of $d$-folding maps on the one-simplex. 
As an example, we computed the first six minimal polynomial folding maps (see Table 1). It is easy to notice that 
these polynomials are the normalized Chebyshev polynomials of the first kind
$P_{1,d}(x)=\frac{1-\cos(d \arccos (1-2x))}{2}$. This is indeed expected from the definition of the Chebyshev
polynomials, which employs the piecewise linear folding function \cite{AR64}. 

\subsubsection{Folding maps of the two-simplex}

From the theory of affine Weyl groups, we know that every tessellation of the two-dimensional plane
that is derived from the reflections of an isosceles triangle is isomorphic to one of two possibilities (see Fig. 4.). 
We denote each possibility by the name of its corresponding  affine Weyl group ($A_2$ and $B_2$).
Of these, only the folding maps derived from the $B_2$ tessellation admit polynomial representatives.
This is because the folding maps derived from the $A_2$ tessellation map every edge in the 
boundary of the triangle to more than one different edge, and therefore the pre-image of an edge cannot be an
algebraic set.

\begin{figure}[hpbt]
\begin{center}
\vspace{0.5cm}
\begin{tabular}{cc}
    \includegraphics[width=7.0cm, angle=90]{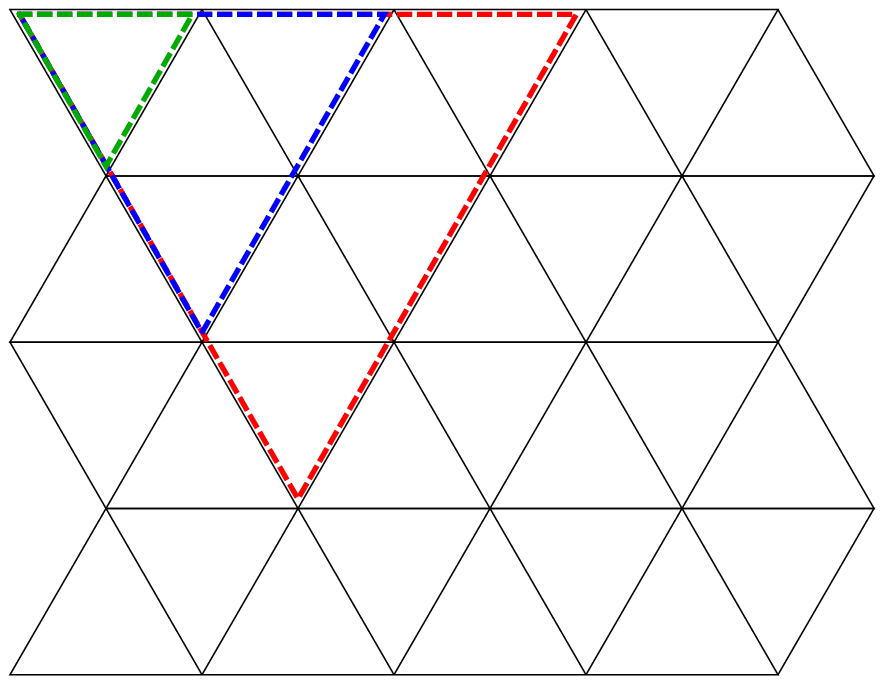}&
    \includegraphics[width=7.0cm, angle=90]{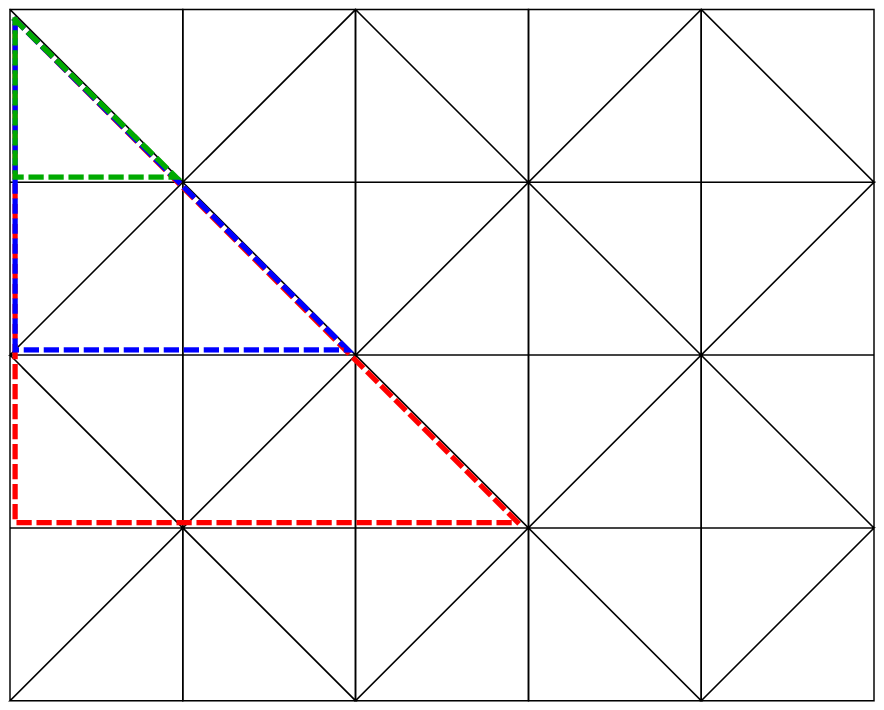}\\
\end{tabular}
\caption{\label{figlie} Two tessellations of the real plane by $2$-simplices. 
The tessellation on the left is associated with the Weyl group $A_2$ and the equilateral triangle.
The right picture represents a tessellation by isosceles right triangles, which is associated with the Weyl group $B_2$.
Both tessellations give rise to corresponding piecewise linear folding maps on the $2$-simplex. 
The dashed green triangles associated with the $g=1$ dilation, give rise to the $1$-fold; the
dashed blue triangles correspond to $g=2$ dilations associated with $4$-folds; the red $g=3$
dilations give rise to $9$-folds. Of these folds only the ones
associated with the isosceles right triangles admit fc-homotopic polynomial maps.
}
\end{center}
\end{figure}

Given the fc-homotopy class associated with the $B_2$ two-fold, using Theorem 2, and assuming that the minimal 
polynomial folding map has degree $k=2$, the defining polynomials have to be of the form
\begin{equation}
\label{twofold}
\left( \begin{array}{c}
P_1(x,y)\\
P_2(x,y)\\
1-P_1(x,y)-P_2(x,y)
\end{array} \right) =
\left( \begin{array}{c}
A(y-Bx)^2 \\
(1-x-y)(C+Dx+Ey)\\
Fxy
\end{array} \right).
\end{equation}
Here, $A,B,\ldots, F$, are unknown parameters that satisfy 
$A>0$, $B>0$, $F>0$ and $C+Dx+Ey>0$ for all $(x,y)\in \Delta^2$.
Now, we impose Eq. \eqref{defeq} on this family of polynomials, which gives rise to the equation
$$
Fxy=1-A(y-Bx)^2-(1-x-y)(C+Dx+Ey).
$$
Expanding this equation in the standard basis of polynomials of degree two with two variables, we get $0=$
$$
(1-C)+(C-D)x+(C-E)y+(E+D+2AB-F)xy+(D-AB^2)x^2+(E-A)y^2.
$$
Therefore, finding the parameters $(A,B,\dots, F)$ that satisfy this equation for all $(x,y)\in \Delta^2$ 
amounts to solving the following system of polynomial equations:
\begin{equation}
\begin{array}{c}
1-C=0, \\
C-D=0, \\
C-E=0,\\
E+D+2AB-F=0, \\
D-AB^2=0, \\
E-A=0.
\end{array}
\end{equation}
From the first, second and third equations it follows that $C=D=E=1$, which yields the strictly positive
linear form $1+x+y\in \mathcal{K}_{+}(\Delta^n,1)$. From the last equation we find that $A=1$, and therefore $F$ and $B$
have to obey the fourth and fifth equations with $A=C=D=E=1$. This implies that
$2=F-2B$ and $1-B^2=0$. As $B$ has to be positive in order to yield a folding map compatible with
the fc-homotopy class, $B$ has to be $1$ and $F$ has to be $4$. This shows that there is a unique polynomial
two-fold of degree two in this fc-homotopy class. We can compactly write this map as
\begin{equation}
\label{twofoldf}
\left( \begin{array}{c}
P_1(x,y)\\
P_2(x,y)\\
1-P_1(x,y)-P_2(x,y)
\end{array} \right) =
\left( \begin{array}{c}
(y-x)^2 \\
(1-x-y)(1+x+y)\\
4xy
\end{array} \right).
\end{equation}

\begin{table}[ht]
\begin{center}
\begin{tabular}{|c|}
\hline\hline
$f_d=$ $\{ P_1(x,y)$, $P_2(x,y)$, $P_3(x,y)=1-P_1(x,y)-P_2(x,y) \}$\\ [0.5ex]
\hline\hline
$f_1=$ $\{x$, $y$,  $1-x-y\}$ \\
\hline
$f_2=$ $\{(x-y)^2$, $(1-x-y)(1+x+y)$,  $4xy\}$ \\ 
\hline
$f_4=$ $\{(1-2x^2-2y^2)^2$, $8xy(1-2xy)$,  \\
$ $ $4(1-x-y)(1+x+y)(x-y)^2\}$ \\
\hline
$f_8=$ $\{\left( 1-4x^2+4x^4-8xy-4y^2+24x^2 y^2 + 4y^2 \right)^2$,\\
$ $ $8(1-x-y)(1+x+y)(1-2x^2+2x^4+4xy-2y^2-4x^2 y^2+2y^4)(x-y)^2$, \\
$ $ $32xy(1-2xy)(1-2x^2-2y^2)^2\}$ \\
\hline
$f_9=$ $\{x(1-y)(2-9x+24x^2-16x^3+9y-24y^2+16y^3)(3-4x)^2(1-4y)^2$,\\
$ $ $y(1-x)(2-9y+24y^2-16y^3+9x-24x^2+16x^3)(3-4y)^2(1-4x)^2$, \\ 
$ $ $(1-x-y)^2\left( 1-8x+16x^2-8y-16xy+16y^2 \right)^2\}$ \\
[1ex]
\hline
\end{tabular}
\caption{\label{tableB} Polynomial factorization of five different folding maps from the 
$2$-simplex to itself.}
\end{center}
\end{table}

One can repeat this construction for higher degree folding maps (see Table 2 and Fig. 5).
The minimal polynomial folding maps of degree four and eight are unique 
and they can be expressed as a composition of the two-folding map; i.e. $f_4=f_2\circ f_2$ and $f_8=f_2\circ f_2 \circ f_2$.
In the case of the nine-folding map we only explored the space of maps corresponding to a single partition of the degrees:
$\{ (2,4),\, (2,4),\, (2,4)\}$. In this partition the interior pre-images of the boundary correspond to quadratic curves (see last row in Table 2 and Fig. 5).
Although we found a unique minimal folding map, other different partitions of the degrees in the defining polynomials
could give rise to additional minimal maps. The fact that there may exist more than one equivalent minimal folding map
can be exemplified by two folding maps of degree eighteen: $f_9\circ f_2$ and $f_2\circ f_9$, which are at the same time distinct and fc-homotopically 
equivalent maps.

\begin{figure}[hpbt]
\begin{center}
\vspace{0.5cm}
\begin{tabular}{cc}
    \includegraphics[width=4.0cm]{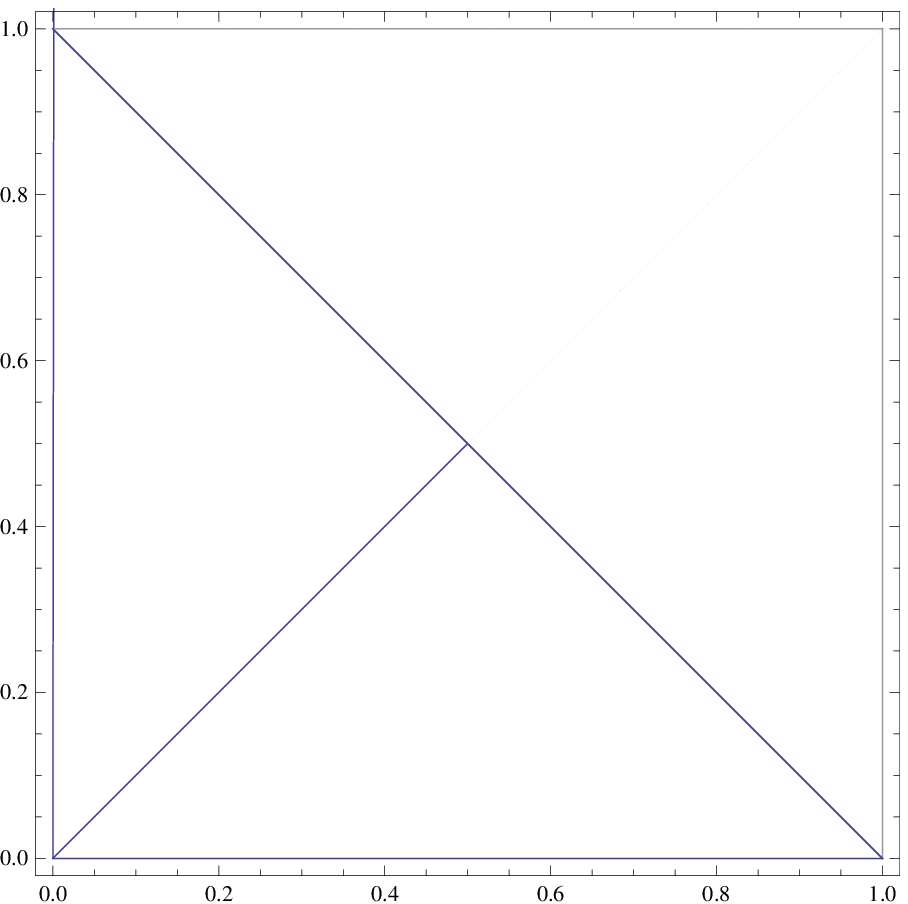}&
    \includegraphics[width=4.0cm]{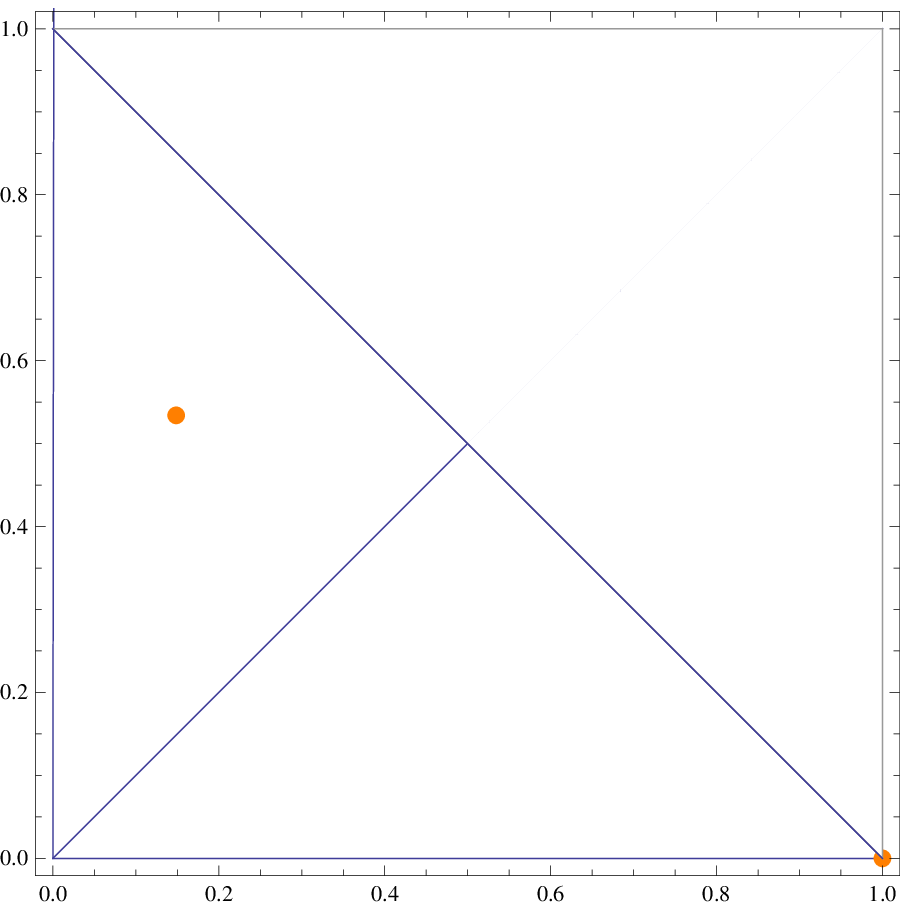}\\
    \includegraphics[width=4.0cm]{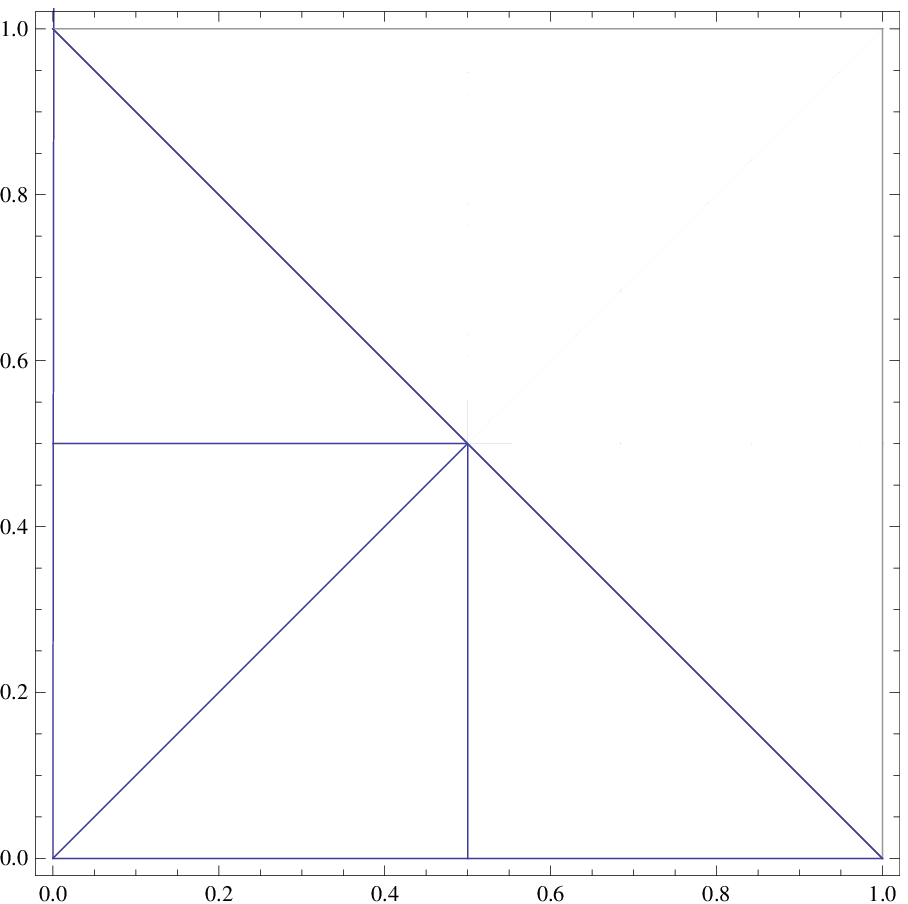}&
    \includegraphics[width=4.0cm]{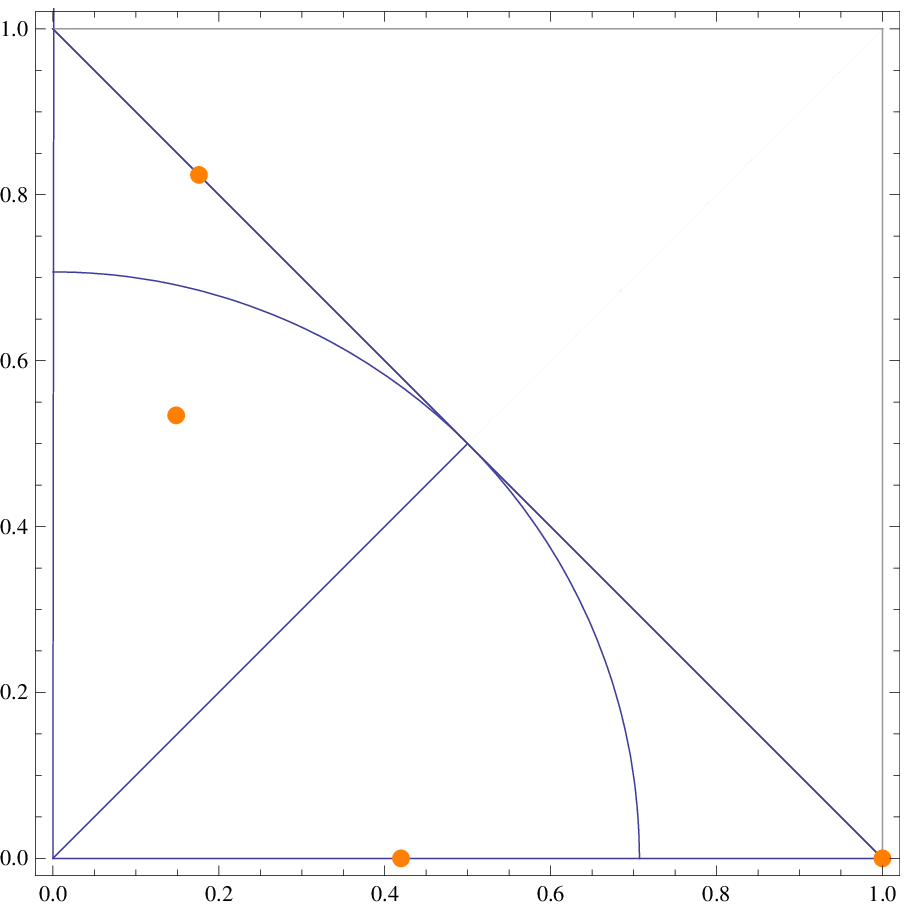}\\
    \includegraphics[width=4.0cm]{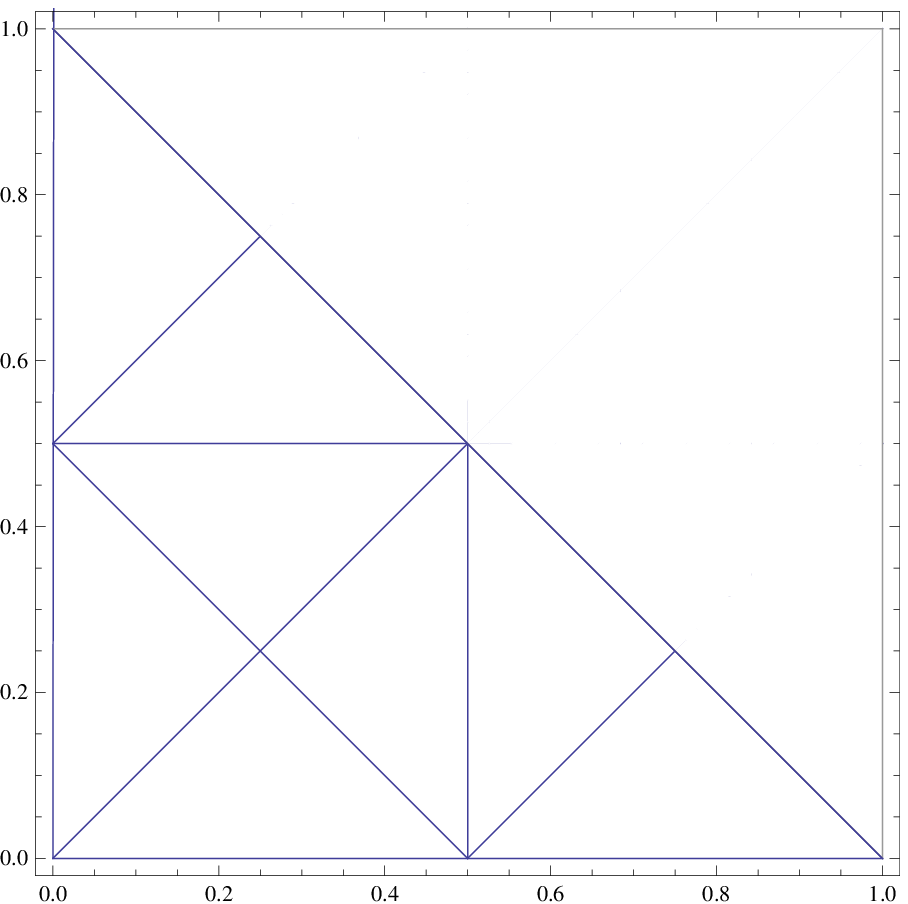}&
    \includegraphics[width=4.0cm]{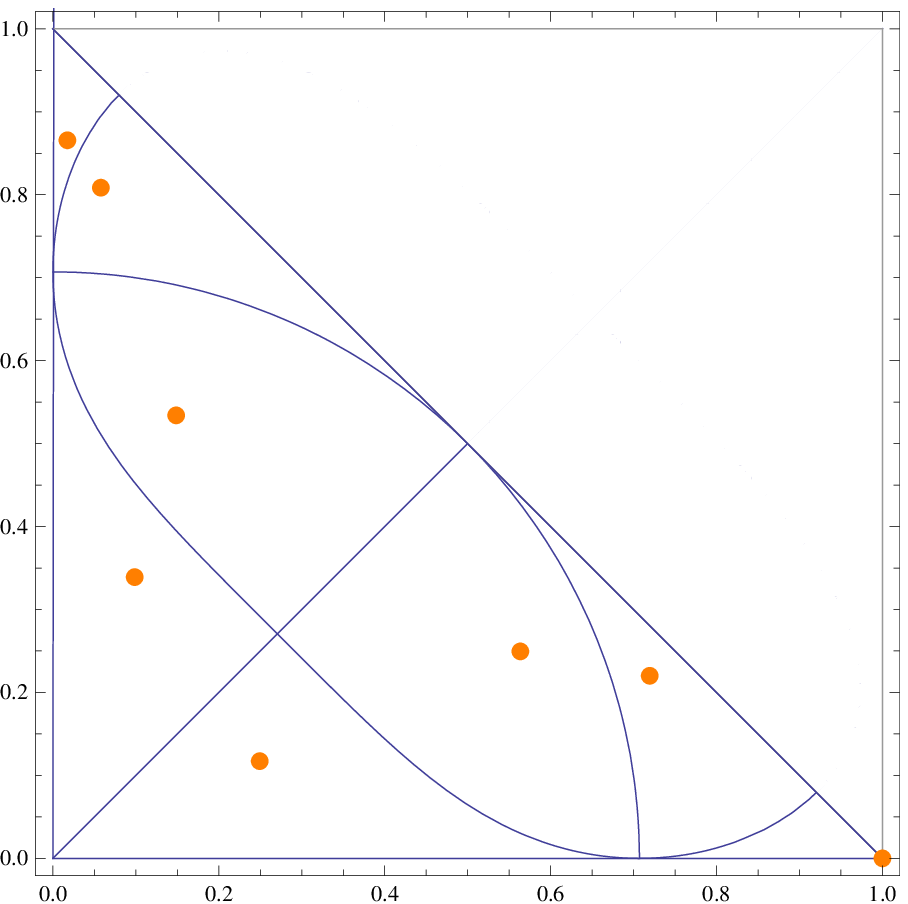}\\
    \includegraphics[width=4.0cm]{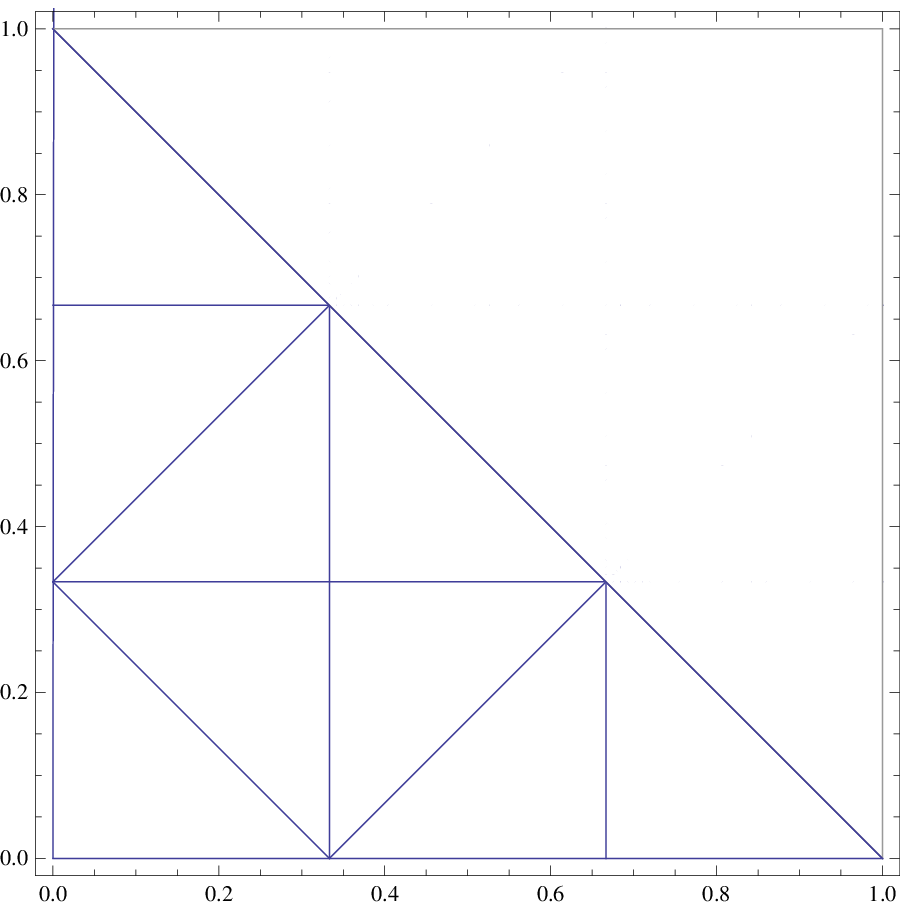}&
    \includegraphics[width=4.0cm]{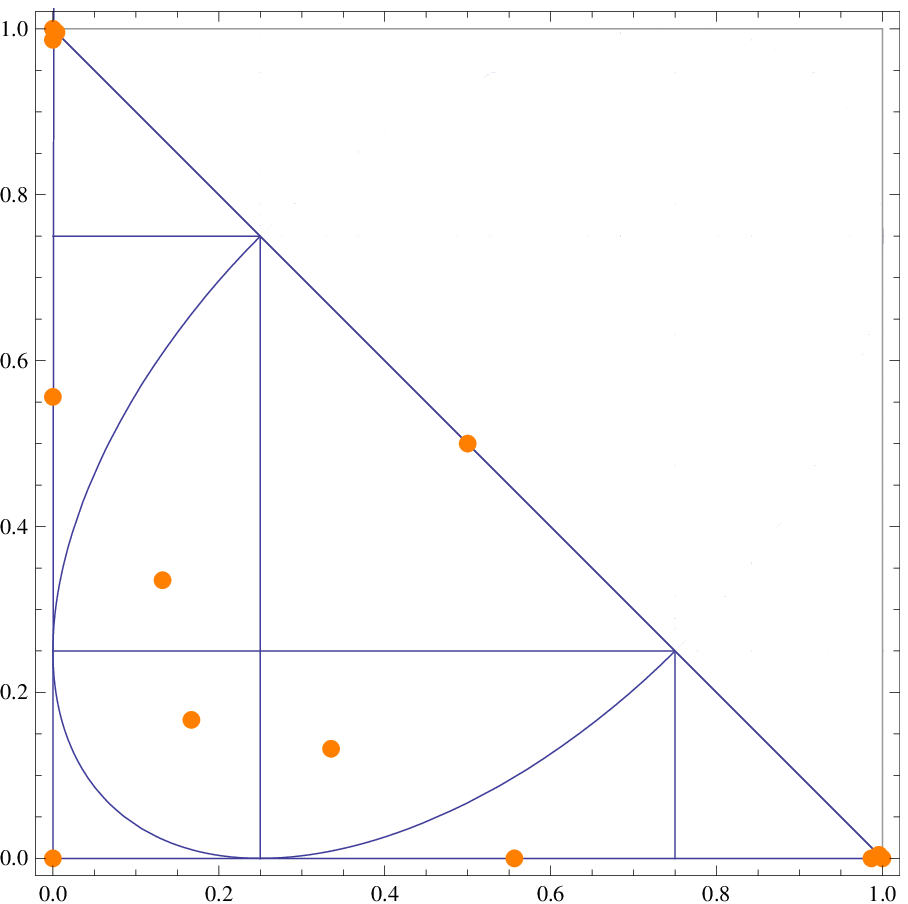}\\
\end{tabular}
\caption{\label{foldinglines} Several folding maps on the $2$-simplex. The lines within each triangle denote
the different preimages of the boundary of the simplex under a given folding map.
Each row consists of an fc-homotopy class of folding maps 
represented by a piecewise linear map (see left)
and a minimal polynomial folding map (see right and Table 2).
The orange dots denote the fixed points of the minimal polynomial folding maps.
}
\end{center}
\end{figure}

\subsubsection{Dynamical properties of some polynomial folding maps}

The set of Chebyshev maps on the one-simplex has very rich dynamical properties \cite{AR64}. 
First, one can show using simple trigonometry that every Chebyshev map of degree $d$ preserves the measure
$$
\mu(B) =\frac{1}{\pi}\int_B \frac{dx}{\sqrt{x(1-x)}},\quad\, B\subseteq [0,1].
$$
In addition to this, the corresponding dynamical systems are strongly mixing with respect to $\mu$; they form a commutative 
semigroup under the composition of maps ($f_d\circ f_e=f_e\circ f_d=f_{d\times e}$); and they are orthogonal as polynomial functions 
under the inner product $L^2([0,1],\mu)$ defined by the invariant measure $\mu$ \cite{AR64}. Some of these properties are also present in 
deformations of the Chebyshev maps that are not folding maps.
For instance, it is well known that many ergodic properties of the logistic map $f_2$ are preserved by  
deformations of the type $f_2^{(\lambda)}=\lambda x(1-x)$, where $\lambda$ belongs to a subset of full measure
in the interval $\sim (3.57, 4]$ (see \cite{May76}).

Despite these results in dimension one, our knowledge regarding polynomial folding maps on higher dimensional simplices is non-existent.
Because of this, we performed several numerical experiments with the polynomial folding maps of degree two
and nine that we derived above (Table 2). In the case of the polynomial two-fold, we determined the fixed points
of the first three iterations of the map (i.e., $f_2$, $f_4$ and $f_8$). We found that the spectrum of the Jacobian matrix
associated with these fixed points lies outside the unit disk in the complex plane, which is suggestive of chaotic 
behavior. We iterated $10^4$ times the two-folding map for different random initial conditions, and found that the dynamics
never converged to an attracting periodic orbit. We repeated this numerical analysis for thousands of maps in the neighborhood
of $f_2$ in $\mathsf{PMaps}(\Delta^2,2)$, and found that the closest maps to $f_2$ exhibited either ergodic behavior or had
an attracting periodic orbit of size larger than $10^4$ iterations (see Fig. 6). 
In order to perform these experiments, we endowed $\mathsf{PMaps}(\Delta^2,2)$ with an $L^2$ metric, and computed 
deformations of the two-folding map by constructing explicitly
a polytopic approximation of the interior of $\mathsf{PMaps}(\Delta^2,2)$ (see Appendix B for details).

\begin{figure}[hpbt]
\begin{center}
\vspace{0.5cm}
\begin{tabular}{c}
    \includegraphics[width=11.0cm]{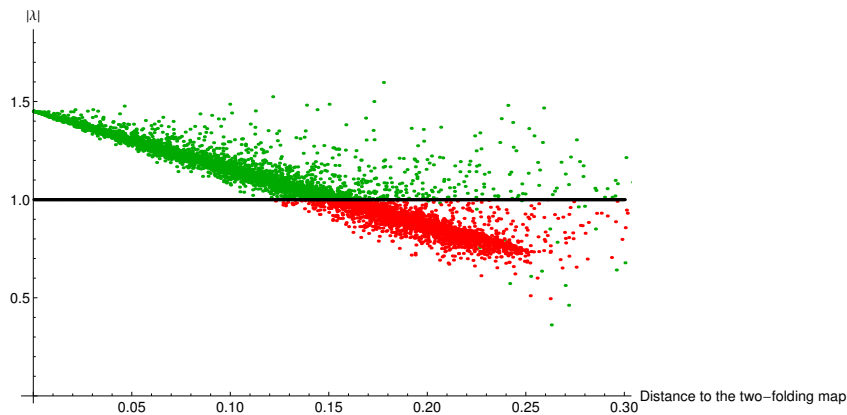}\\
\end{tabular}
\caption{\label{figtwofold} Dynamical properties of deformations of the polynomial two-folding map
on the $2$-simplex (see Appendix B). The horizontal axis denotes the $L^2$ distance to the two-folding map. 
The vertical axis denotes the smallest absolute value of the eigenvalues of the Jacobian matrix associated with the fixed points of the map. 
The green dots denote maps that are either ergodic or their dynamics converge to a periodic orbit
of size larger than $10^4$ iterations. The red dots denote maps whose dynamics converge either to a finite set of absorbing fixed points
or to a periodic orbit of size smaller than $10^4$ iterations.
}
\end{center}
\end{figure}

In the case of the nine-folding polynomial map we found thirteen fixed points. The Jacobian matrices associated with
eleven of these fixed points have their spectrum located in the exterior of the unit disk, while in the remaining fixed points, $(1,0)$ and $(0,1)$,
the Jacobian matrix has eigenvalues in the interior of the unit disk. These vertices are attracting fixed points, and each one is located
next to two repelling fixed points (within a Euclidean distance of less than $1/100$). We computed the dynamics of $10,000$ random initial
conditions close to the origin $(0,0)$ and computed their fixation times (see Fig. 7). As expected, every scenario converged either to
the vertex $(0,1)$ or $(1,0)$, and the distribution of fixation times closely follows a log-normal distribution with mean $2.557$
and standard deviation $0.482$.

\begin{figure}[hpbt]
\begin{center}
\vspace{0.5cm}
\begin{tabular}{c}
    \includegraphics[width=10.0cm]{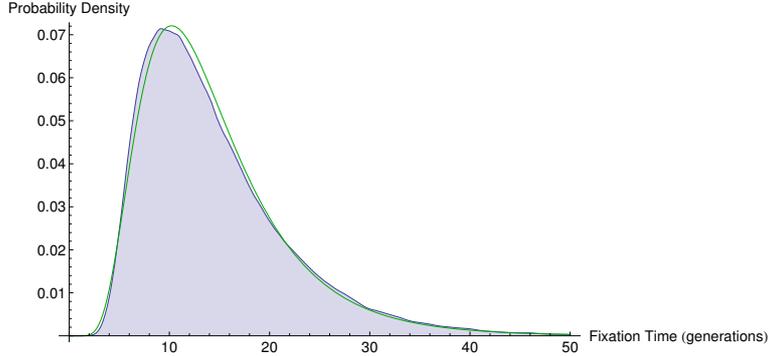}\\
\end{tabular}
\caption{\label{figninefold} Distribution of fixation times associated with the dynamics of a minimal polynomial nine-fold
(see Table 2). Here, time is measured as the number of iterations.
The vertex points $(1,0)$ and $(0,1)$ of the triangle (see Fig. 5) are the only stable fixed points of the map.
The distribution represents the fixation times of
$10^4$ uniformly distributed random points in the square $0<x<1/100$, $0<y<1/100$. The shape of the resulting
distribution did not change appreciably when smaller squares were chosen. 
A fit of the log-normal distribution ($\hat{\mu}=2.557, \hat{\sigma}=0.482$) to the simulated data (green curve) is shown superposed.
}
\end{center}
\end{figure}

\section{Conclusion}

Discrete-generation multi-genic models with non-trivial 
linkage maps are difficult to study because of the nonlinearity associated with the dynamical equations. 
For instance, simple neutral mutation-recombination models with no genetic drift 
require quadratic maps from the simplex to itself to describe the corresponding population dynamics. 
More generally, the addition of weak selection effects combined with time-dependent intensities of selection, mutation
and recombination can give rise to very complex dynamical systems also described by the iteration of
polynomial maps from the simplex to itself  (e.g. see \cite{Lyub92, Gro07, Baa01, NHB99}).
Although these maps introduce a plethora of technical problems that are not present
in the case of finite Markov chains, some particular solutions can be found in biologically
interesting regimes, such as the loose-linkage limit \cite{NHB99} and the tight-linkage low-mutation limit
\cite{Baa01}.

In this paper, rather than stressing particular models, we have described some general properties of the space
of stochastic polynomial maps on the simplex. Our motivation here is twofold. On the one hand we are interested
in developing efficient parametrizations of the space of polynomial maps of bounded degree on the simplex, 
with the goal of providing a large class of diffeomorphisms on the simplex that are useful in practical 
applications (e.g. defining adaptive meshes on the simplex). On the other hand, a global description of
the space of polynomial maps allows us to systematically study properties of the associated dynamical systems.
For instance, we have characterized an important class of maps in the boundary of $\mathsf{PMaps}(\Delta^{n},k)$,
the minimal polynomial folding maps. 
We have shown how the dynamics of some of these folding maps and their deformations in $\mathsf{PMaps}(\Delta^{n},k)$ 
can exhibit ergodic and mixing properties, similar to the logistic map in dimension one \cite{May76}.
Still, it would be interesting to explore whether there exist more efficient algorithms to construct 
these folding maps. The construction that we have provided in this paper requires 
significant computational resources. However, given that the one-dimensional 
minimal polynomial folding maps 
can be derived easily using classical methods, it remains an open 
question whether the higher dimensional folding maps can be derived using simpler methods than the ones that we use here.

\section*{Acknowledgements}
We thank Ethan Akin, Ben Greenbaum, Julien Keller, Stan Leibler, Arnold Levine, Rami Pugatch, Tiberiu Tesileanu and Tsvi Tlusty for helpful comments. 
S.L. is the Addie and Harold Broitman Member
in the Simons Center for Systems Biology at the Institute for Advanced Study, Princeton, NJ.

\newpage
\appendix
\section{Some Elementary Lemmas}

\begin{lemma}
Every one-to-one and onto stochastic matrix in $\mathsf{PMaps}(\Delta^{n},1)$ is an $n+1$ by $n+1$ permutation matrix.
\end{lemma}
\emph{Proof.} Let $M$ be a one-to-one and onto stochastic matrix, represented by 
an $n+1$ by $n+1$ invertible matrix. Let the standard simplex $\Delta^n$ be defined 
as the convex hull of the standard basis elements in $\mathbb{R}^{n+1}$, 
i.e. $e_1=(1,0,0,\cdots)$, $e_2=(0,1,0,\cdots)$, $\ldots\, e_{n+1}=(0,0,0,\cdots,1)$.
Therefore, the image of $\Delta^n$ under $M$ in $\mathbb{R}^{n+1}$ is the convex hull of the vertices
$\{ Me_i \}_{i=1}^{n+1}$. As $M$ is a stochastic matrix, the image of the standard simplex is a subset of the standard simplex.
Furthermore, as $M$ is onto, the image simplex is indeed the standard simplex, i.e. the convex hull of the standard basis elements. 
This implies that $Me_i=e_j$ for every $i\in [1,n+1]$ and
 $j\in [1,n+1]$. As $M$ is invertible, the image of the standard basis elements under $M$ spans all the standard basis elements, 
and therefore $M$ has to be a permutation matrix.$\square$
\vspace{0.75cm}

\begin{lemma}
Let $f\colon\, \Delta^n\to\Delta^n$ be a folding map. The image of the boundary of the simplex under $f$
is contained in the boundary, $f(\partial\Delta^n)\subset \partial\Delta^n$.
\end{lemma}
\emph{Proof.} It follows from the definition of folding map that $f$ is continuous, onto and $d$-to-$1$, such that the number
of preimages is $d$ everywhere in the interior of $\Delta^n$.
Let us assume that there exists a point $x\in \partial\Delta^n$ whose image under $f$ is in the interior of the simplex. Then, one can 
construct a small open ball $B_{f(x)}$ also in the interior of the simplex that contains $f(x)$. As $f(x)$ is the image
of a point in the boundary of $\Delta^n$ we can divide $B_{f(x)}$ in three subsets. 
First, one subset is an open disk of dimension $n-1$ that is the image under $f$ of a subset in the boundary
of $\Delta^n$ that contains $x$. This $n-1$ dimensional disk, that we denote by $D_{f(x)}$, divides the ball in two open half-balls
$B_{f(x)}^+$ and $B_{f(x)}^-$, such that $B_{f(x)}=B_{f(x)}^+\cup D_{f(x)} \cup B_{f(x)}^-$.
By the continuity of $f$, one of the half-balls $B_{f(x)}^+$ has a preimage open subset in $\Delta^n$ whose boundary contains
$x\in\Delta^n\subset\mathbb{R}^n$. However, the other half-ball $B_{f(x)}^-$ does not have a preimage 
in $\Delta^n$ whose boundary contains $x$. It follows that the number of preimages of $x$ is larger than the number
of preimages of points in $B_{f(x)}^-$. This is a contradiction because $f$ is a folding map,
and by definition the number of preimages is the same everywhere in the interior of the simplex. $\square$
\vspace{0.75cm}

\section{Deformations of maps}

In order to know whether the dynamical properties of a given folding map $f_\ast$ are preserved
in its neighborhood in $\mathsf{PMaps}(\Delta^{n},k)$, we need a systematic method to generate
deformations of $f_\ast$ and determine their properties. 
First, we endow $\mathsf{PMaps}(\Delta^{n},k)$ with the $L^2$ distance associated with 
the Euclidean volume form on the simplex, i.e.
$$
\Vert f_\ast-g \Vert_{L^2} = \sqrt{ \int_{\Delta^n}\left[\sum_{i=1}^n (P_i(x)-Q_i(x))^2 \right] \prod_{i=1}^n dx_i },
$$
where $f_\ast ,\, g\in \mathsf{PMaps}(\Delta^{n},k)$, and $\{ P_i(x) \}_i$ and $\{ Q_i(x) \}_i$ are the sets 
of defining polynomials for $f_\ast$ and $g$ respectively.
Equipped with this metric we can generate deformations $g\in B_{f_\ast,\epsilon}$ of the folding map within a 
distance $\epsilon$,
$$
B_{f,\epsilon} =\{ g\in\mathsf{PMaps}(\Delta^{n},k)\colon \,\, \Vert f_\ast-g \Vert_{L^2} \leq \epsilon\},
$$
by following the steps outlined in Fig. 8. In particular, we first generate random interior points $r\in$ 
$\mathsf{PMaps}(\Delta^{n},k)$ and then take convex combinations of those points with the folding map
to obtain the desired deformations $g=tr+(1-t)f_\ast\in B_{f_\ast,\epsilon}$. The $t$ parameter is chosen
at random from the uniform distribution on the interval 
$(0,t_\epsilon]$, with $t_\epsilon = \epsilon/\Vert f_\ast-r \Vert_{L^2}$.

\begin{figure}[hpbt]
\begin{center}
\vspace{0.5cm}
\begin{tabular}{c}
    \includegraphics[width=9.0cm]{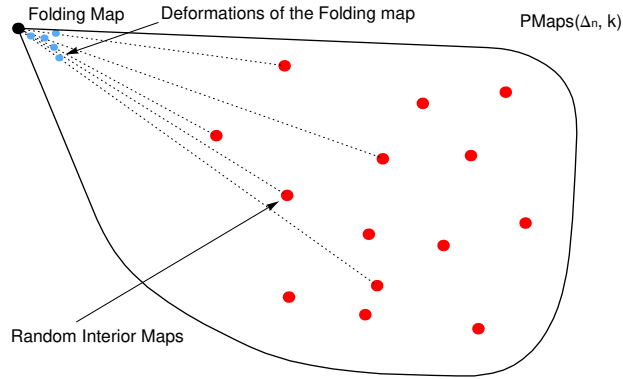}\\
\end{tabular}
\caption{\label{figdefmaps} Schematic representation of our method to construct deformations of a folding map. First, we
construct random maps in the interior of the space of polynomial maps (red dots). Then, we use the convexity of 
$\mathsf{PMaps}(\Delta^{n},k)$ to construct random convex combinations of interior points and the folding map such 
that the resulting deformations are ``close'' to the folding map (blue dots).
}
\end{center}
\end{figure}

As we showed in the proof of Theorem 1, $\mathsf{PMaps}(\Delta^{n},k)$ is a bounded convex subset
of the cone $\times_{i=1}^n \mathcal{K}(\Delta^n,k)\subset \oplus_{i=1}^n \Pi^{n}_k$. 
Hence, in order to parametrize the interior of the space of polynomial maps we need to 
first parametrize the convex cone of strictly positive polynomials $\mathcal{K}_{+}(\Delta^n,k)$, whose closure is the cone
of non-negative polynomials $\mathcal{K}(\Delta^n,k)$.
This is difficult to attain in practice, because $\mathcal{K}_{+}(\Delta^n,k)$ is not finitely generated for values of $k$ larger than one. 
However, using P\'olya's criterion one can construct an infinite sequence of finitely generated cones $\{ \mathcal{K}_{N}(\Delta^n,k) \}_{N=1}^\infty$
that satisfy
\begin{equation}
\mathcal{K}_{N}(\Delta^n,k) \subseteq \mathcal{K}_{N+1}(\Delta^n,k) \subseteq \mathcal{K}_{N+2}(\Delta^n,k) \subseteq  \cdots \to \mathcal{K}_{+}(\Delta^n,k),
\label{eqnflag}
\end{equation}
and converge to $\mathcal{K}_{+}(\Delta^n,k)$ in the limit of large $N$. Therefore, in the numerical applications that concern us here, we use the cone
$\mathcal{K}_{N}(\Delta^n,k)$ with large $N$ as an approximation of $\mathcal{K}_{+}(\Delta^n,k)$. 
This can be done more precisely as follows. First, we consider the standard parametrization of the vector space
of polynomials $\Pi^{n}_k$, in which every polynomial $P(x)$ is expressed as a linear combination of monomials
\begin{equation}
\label{expansionP}
P(x)=\sum_{\vert\alpha\vert\leq k} c_{\alpha_1\alpha_2\cdots\alpha_n}x_1^{\alpha_1}x_2^{\alpha_2}\cdots x_n^{\alpha_n}.
\end{equation}
Then we construct $P_H(x)$, the homogenization of $P(x)$, which is obtained by multiplying every term of degree $\vert\alpha\vert<k$ in Eq. \eqref{expansionP}
by the factor $\left( \sum_{i=1}^{n+1} x_i \right)^{k-\vert\alpha\vert}$. Finally, we multiply $P_H(x)$ by $( \sum_{i=1}^{n+1} x_i)^N$ and
expand the resulting polynomial 
$$
P_H(x)\times\left( \sum_{i=1}^{n+1} x_i \right)^N
$$ 
in the basis of homogenous monomials of degree $N+k$ in $n+1$ variables:
$$
\left\{ x_1^{\beta_1}x_2^{\beta_2}\cdots x_n^{\beta_n}x_{n+1}^{\beta_{n+1}} \right\}_{\vert \beta\vert=N+k}.
$$
The reader can notice that each term in this expansion corresponds to a linear combination of the $c_{\alpha}$-coefficients 
that appear in Eq. \eqref{expansionP}. We then apply P\'olya's theorem, which requires each coefficient of the expansion to be 
positive \cite{HLP52}. As the dimension of the space of homogeneous polynomials of degree $N+k$ in $n+1$ variables
is $(N+k+n)!/((N+k)!n!)$, imposing that each coefficient in the expansion has to be positive
gives rise to an equal number of inequalities in $\Pi^{n}_k$. 
We define $\mathcal{K}_{N}(\Delta^n,k)$ to be the semi-algebraic subset that satisfies these inequalities.
Furthermore, in order to parametrize this convex cone one can determine its generators by means of a
vertex enumeration algorithm (e.g. \cite{Avi00}) such that any point in the cone can be expressed as a linear combination 
of the generators with positive coefficients.

In order to learn how one can implement this algorithm, it is useful to consider a particular example.
Here, we consider deformations of the polynomial two-fold in $\mathsf{PMaps}(\Delta^{2},2)$. In this case, the space of quadrics has dimension six.
To determine the plot shown in Fig. 6, we used the value $N=8$ to construct the cone  $\mathcal{K}_{8}(\Delta^2,2)$.
The application of P\'olya's criterion gave rise to $66$ inequalities in $\mathbb{R}^6\simeq \Pi_2^2$, which 
define a cone with $900$ generators 
(we used the reverse search algorithm $\mathsf{lrs}$ to determine the generators of the cone \cite{Avi00}). We scaled each generator $e_i$ 
by a positive real number $\lambda_i$, such that the maximum value of the corresponding polynomial on the simplex was one. 
Given this basis of scaled generators, we used the Dirichlet distribution with $901$ concentration parameters
$$
\alpha_1=\alpha_2=\cdots=\alpha_{901}=1/1000,
$$
to generate random convex combinations of the scaled generators. These combinations correspond to rays in 
$\mathcal{K}_{8}(\Delta^2,2)$. Furthermore, by choosing a random number sampled from the uniform distribution on 
$(0,1)$ we generate random quadrics in the segment of the ray defined by the origin of the cone and 
the random convex combination of the scaled generators. We applied this algorithm twice to generate
pairs of random quadrics $\{ R^\ast_1(x),\, R^\ast_2(x) \}$. Then, 
we computed the maximum of the function $R^\ast_1(x)+R^\ast_2(x)$ on $\Delta^2$, that
we denote by $S_{max}$. It follows that sampling a random number $t^\ast$ from the uniform distribution
on $[0,1/S_{max}]$, defines an interior map $r\in \mathsf{PMaps}(\Delta^{2},2)$ with defining polynomials
$$
R_1(x)=t^\ast R^\ast_1(x),\quad R_2(x)=t^\ast R^\ast_2(x).
$$
Finally, the construction of a deformation $g$ of the two-folding map 
$f_\ast$ requires an additional convex combination of $f_\ast$ with the map $r$ (see Fig. 8.).

\end{document}